\documentclass{amsart}
\usepackage{amssymb,amscd}

\numberwithin{equation}{section}

\newtheorem{prop}{Proposition}
\newtheorem{theorem}[prop]{Theorem}
\newtheorem{corollary}[prop]{Corollary}
\newtheorem{lemma}[prop]{Lemma}

\theoremstyle{definition}
\newtheorem{definition}[prop]{Definition}
\newtheorem{example}[prop]{Example}
\newtheorem{remark}[prop]{Remark}

\numberwithin{prop}{section}

\newfont{\germ}{eufm10}

\newcommand\B{{\mathcal B}}
\newcommand\et[1]{\tilde{e}_{#1}}
\newcommand\ft[1]{\tilde{f}_{#1}}
\newcommand\geh{\mbox{\germ g}}

\newcommand\ot{\otimes}

\newcommand\veps{\varepsilon}
\newcommand\vphi{\varphi}

\newcommand\xb{\overline{x}}
\newcommand\Z{\mathbb{Z}}

\newcommand\yb{\overline{y}}

\begin{document}

\title{Geometric Crystal and Tropical $R$ For 
$D^{(1)}_n$}

\author{A. Kuniba}
\address{Institute of Physics, University of Tokyo, Tokyo 153-8902, Japan}
\email{atsuo@gokutan.c.u-tokyo.ac.jp}

\author{M. Okado}
\address{Department of Informatics and Mathematical Science,
    Graduate School of Engineering Science,
Osaka University, Osaka 560-8531, Japan}
\email{okado@sigmath.es.osaka-u.ac.jp}

\author{T. Takagi}
\address{Department of Applied Physics, National Defense Academy,
Kanagawa 239-8686, Japan}
\email{takagi@nda.ac.jp}

\author{Y. Yamada}
\address{Department of Mathematics, Faculty of Science,
Kobe University, Hyogo 657-8501, Japan}
\email{yamaday@math.kobe-u.ac.jp}


\begin{abstract}
We construct a geometric crystal 
for the affine Lie algebra $D^{(1)}_n$ 
in the sense of Berenstein and Kazhdan.
Based on a matrix realization including a spectral parameter, 
we prove uniqueness and explicit form of the 
tropical $R$, the birational
map that intertwines products of the geometric crystals.
The tropical $R$ commutes with geometric Kashiwara operators and 
satisfies the Yang-Baxter equation.
It is subtraction-free and yields a piecewise
linear formula of the combinatorial $R$ for crystals 
upon ultradiscretization.
\end{abstract}

\maketitle

\section{Introduction}
In this paper we construct a  
geometric crystal for $D^{(1)}_n$ corresponding to 
the $D^{(1)}_n$-crystal in \cite{KKM}, and 
study the associated tropical $R$.
Geometric crystal is a notion introduced by 
Berenstein and Kazhdan in \cite{BK} as an 
algebro-geometric analogue of the crystal theory \cite{K1,K2}.
The latter is a theory of quantum groups at $q=0$, which 
involves combinatorial operations described by 
various piecewise linear functions.
In the geometric setting in \cite{BK}, 
such a structure corresponds to birational maps 
described by subtraction-free functions.

Let us illustrate these features in the case of $A^{(1)}_{n-1}$
based on the constructions in \cite{BK} and \cite{KNY}.
Elements $x$ of the crystal $B_l$ for the $l$-symmetric tensor representation
is specified by the coordinates 
$x=(x_1,\ldots, x_n) \in \Z_{\ge 0}^n$ with $x_1 + \cdots + x_n = l$,
where all the indices are considered to be in $\Z/n\Z$ \cite{KKM}.
The action of the Kashiwara operator 
${\tilde e}_i$ up to 2-fold tensor product is given, unless they vanish, by
\begin{align*}
& {\tilde e}^c_i(x) = (\ldots, x_{i-1}, x_i+c, x_{i+1}-c, x_{i+2},\ldots),\\
& {\tilde e}_i^c(x \otimes y) 
= {\tilde e}^{c_1}_i(x) \otimes {\tilde e}_i^{c_2}(y),\\
&c_1 = \max(x_i+c,y_{i+1}) - \max(x_i,y_{i+1}), \\
&c_2 = \max(x_i,y_{i+1}) - \max(x_i,y_{i+1}-c).
\end{align*}
In the geometric crystal, one still has the coordinates 
$x=(x_1,\ldots, x_n) \in \B$ and the corresponding structure looks 
\begin{align*}
&e^c_i(x) = (\ldots, x_{i-1}, cx_i, c^{-1}x_{i+1}, x_{i+2},\ldots),\\
&e_i^c(x, y) 
= (e^{c_1}_i(x), e_i^{c_2}(y)),\\
&c_1 = \frac{cx_i+y_{i+1}}{x_i+y_{i+1}}, \;\;  
 c_2 = \frac{x_i+y_{i+1}}{x_i+c^{-1}y_{i+1}}.
\end{align*}
We call it the geometric Kashiwara operator.
Note that the latter $c_1, c_2$ contain no minus sign, i.e., 
they are subtraction-free.
The former $c_1, c_2$  are piecewise linear 
and obtained from the latter by replacing 
$+,\times, /$ with $\max, +, -$, respectively.
The procedure, which we call {\it ultradiscretization} in this paper,  
is well defined since the two sides of $p(q+r) = pq+pr$ have the coincident image
$p+\max(q,r) = \max(p+q, p+r)$.  
The geometric crystal $\B$ is actually 
associated with the {\it coherent family}  
$\{ B_l \}_{l\ge 1}$ \cite{KKM} rather than the single crystal $B_l$.

To each $x \in \B$ we assign the $n$ by $n$ matrix 
\[M(x,z)=(\sum_{i=1}^{n}\frac{1}{x_i}E_{i,i}-
\sum_{i=1}^{n-1}E_{i+1,i}-z E_{1,n} )^{-1}\]
containing a spectral parameter $z$.
Then on the product $M=M(x^{(1)},z)\cdots M(x^{(L)},z)$ 
corresponding to $(x^{(1)},\ldots,x^{(L)}) \in \B^{\times L}$, 
one can realize the geometric Kashiwara operator as the multiplications of 
unipotent matrices:
\[G_i\left(\frac{c-1}{z^{\delta_{i0}}\veps_i(M)} \right)M
G_i\left(\frac{c^{-1}-1}{z^{\delta_{i0}}\vphi_i(M)} \right)=
M(e^{c_1}_i(x^{(1)}),z)\cdots M(e^{c_L}_i(x^{(L)}),z),\]
where
$G_i(a)=E+a E_{i,i+1}$, $(1 \leq i \leq n-1)$,
$G_0(a)=E+a E_{n,1}$ and 
\begin{displaymath}
\veps_i(M)=z^{-\delta_{i0}}
\frac{M_{i+1,i}}{M_{i,i}}\Big\vert_{z=0}, \quad
\vphi_i(M)=z^{-\delta_{i0}}
\frac{M_{i+1,i}}{M_{i+1,i+1}}\Big\vert_{z=0}.
\end{displaymath}
The product $(e^{c_1}_i(x^{(1)}),\ldots,e^{c_L}_i(x^{(L)})) 
\in \B^{\times L}$ corresponds to the element 
${\tilde e}^c_i(x^{(1)} \otimes \cdots \otimes x^{(L)})$ in 
the tensor product crystal
under the ultradiscretization mentioned above.

One of the axioms of the geometric crystal is the relation 
$e^c_i e^{cd}_j e^d_i = e^d_j e^{cd}_i e^c_j$ in the case of
$\langle\alpha^\vee_i,\alpha_j\rangle
=\langle\alpha^\vee_j,\alpha_i\rangle=-1$.
The matrix realization reduces its proof to 
a simple manipulation on the matrices $G_i$'s.

In the crystal theory, the isomorphism of tensor products
$B_l\otimes B_k \simeq B_k \otimes B_l$ is called 
the combinatorial $R$.
It is the quantum $R$ matrix at $q=0$, and plays a fundamental role.
To incorporate it into the geometric setting we consider 
the matrix equation
$M(x,z)M(y,z)=M(x'z)M(y',z)$.
Explicitly, it is a discrete Toda type equation 
$x_iy_i=x_i'y_i'$, $\frac{1}{x_i}+\frac{1}{y_{i+1}}=
\frac{1}{x'_i}+\frac{1}{y'_{i+1}}$.
Given $x=(x_1,\ldots, x_n) \in \B$, 
call  $\ell(x) = x_1 \cdots x_n$ the level.
It turns out that the equation with the level constraints 
$\ell(x)=\ell(y'),\, \ell(y)=\ell(x')$ characterizes the unique solution:
\begin{eqnarray*}
&&x'_i=y_i \frac{P_i(x,y)}{P_{i-1}(x,y)}, \quad
y'_i=x_i \frac{P_{i-1}(x,y)}{P_i(x,y)},\\
&&P_i(x,y)=\sum_{k=1}^{n} \prod_{j=k}^n x_{i+j} \prod_{j=1}^k y_{i+j}.
\end{eqnarray*}
We call the rational map $(x,y) \mapsto (x',y')$ the {\it tropical $R$}, 
and simply denote by $R$.
{}From the characterization 
one can derive the geometric analogues of most important properties 
of the combinatorial $R$, e.g., 
commutativity with geometric Kashiwara operators 
$e^c_iR= Re^c_i$, the inversion relation $R(R(x,y))=(x,y)$ 
showing $R$ is birational, the Yang-Baxter relation and 
so forth. 
Moreover the ultradiscretization of the tropical $R$ yields an 
explicit piecewise linear formula for the combinatorial $R$.

In this paper we show that all the features explained above along  
type $A$ persist also for $D^{(1)}_n$.
What we do is to start from a new geometric crystal 
associated with the coherent family of 
$D^{(1)}_n$-crystals in \cite{KKM}, and develop an 
elementary approach based on explicit calculations. 
It is motivated by several recent works; 
geometrization of canonical bases \cite{BFZ,BK,L2},
tropical combinatorics and 
birational representations of affine Weyl groups \cite{Ki,KNY,NoY,Y}, 
integrable cellular automata associated to crystals 
\cite{HKT1,FOY,HKT2,HKOTY},  discrete soliton equations 
arising as their tropical evolution equation \cite{TTMS, HHIKTT,Y}, 
combinatorial $R$ in terms of the bumping algorithm \cite{HKOT1,HKOT2}, etc.
We hope to report on applications of the present result
to these issues in near future.

The paper is arranged as follows.
In section \ref{sec:B} we begin with the definition and basic facts 
on geometric crystals following \cite{BK}. 
The geometric $D^{(1)}_n$-crystal $\B$ is presented, which 
is parametrized with the $2n-1$ coordinates 
$x=(x_1,\ldots,x_n,\xb_{n-1},\ldots,\xb_1)$.
See section \ref{subsec:piecewise} for the 
precise correspondence to the coordinates in the crystal.
In section \ref{sec:M} we introduce the $2n$ by $2n$ 
matrix $M(x,z)$  that realizes the geometric crystal $\B^{\times L}$.
It is quadratic in $z$ in contrast with the 
essentially linear one for type $A$.
At $z=0$ it admits a factorization (\ref{eq:factorizationofA}) and 
belongs to the lower triangular Borel subgroup of $D_n$ as in
Proposition \ref{pr:MSMS}, whereas at $z= \ell(x)^{-1}$ 
it shrinks to the rank one matrix as in Proposition \ref{prop:meqdpd}.
A key in our approach is to 
make full use of such properties 
controlled by the spectral parameter of affine nature,  
which we have been unable to learn in \cite{BK}.
In section \ref{sec:R} we study the tropical $R$.
{}From Theorem \ref{th:uniqueness} 
its uniqueness is immediate, hence our main task 
is to assure the existence.
We shall show that the explicit birational map 
$(x,y) \rightarrow (x',y')$ 
in Definition \ref{th:tropicalR} is the answer.
By the ultradiscretization we then derive a piecewise linear formula for
the combinatorial $R$ and the associated energy function. 
Appendix \ref{sec:appA} provides an alternative proof of 
the relation $e_i^ce_j^{cd}e_i^d=e_j^de_i^{cd}e_j^c$ on $\B^{\times L}$.
Appendix \ref{sec:proofofaaeqaa} contains 
a proof of $M(x,z)M(y,z)=M(x'z)M(y',z)$ at $z=0$.

We remark that our terminology 
``tropical"  in this paper follows the usage in \cite{Ki,NoY}, 
which is unfortunately the opposite of \cite{BK}.  
The phrase ``ultradiscretization" has been used in many
papers inspired by \cite{TTMS}.

\section*{Acknowledgements}
The authors thank M. Noumi, S. Naito, and 
D. Sagaki for useful discussions.
A.K. thanks Shi-shyr Roan 
for the hospitality at Academia Sinica, 
Taipei during ``Workshop on Solvable Models and Quantum Integrable Systems"
in November 2001, where a part of the work was presented.
M.O. and Y.Y. are partially supported by Grant-in-Aid for 
Scientific Research from Ministry of Education, Culture, 
Sports, Science and Technology of Japan.

\section{Geometric $D^{(1)}_n$-crystal}\label{sec:B}

\subsection{Geometric crystal}\label{subsec:gcrystal}
Let $\geh$ be an affine Lie algebra of simply laced type, i.e., 
$\geh=A^{(1)}_n,D^{(1)}_n$
or $E^{(1)}_{6,7,8}$. Let $I=\{0,1,\ldots,n\}$ be the index set of vertices 
of the Dynkin
diagram of $\geh$. We denote by $\{\alpha_i\mid i\in I\}$ (resp. 
$\alpha^\vee_i\mid i\in I\}$)
the set of simple roots (resp. simple coroots). Note that $\langle 
\alpha^\vee_i,\alpha_j\rangle
=0$ or $-1$.

Following \cite{BK} we introduce a geometric crystal associated to $\geh$. 
Throughout this paper,
{\em variable} means a real variable unless otherwise stated. Let 
$x=(x_1,x_2,\ldots,x_N)$ be
an $N$-tuple of variables. For $i\in I$ let $\veps_i(x),\vphi_i(x)$ be 
rational functions in $x$
and set
\begin{equation} \label{eq:gamma}
\gamma_i(x)=\vphi_i(x)/\veps_i(x).
\end{equation}
For $i\in I$ let $\mathbb{R}\times\mathbb{R}^N\longrightarrow\mathbb{R}^N:
(c,x)\mapsto e_i^c(x)$ be a rational mapping, namely, 
$e_i^c(x)=(x'_1,x'_2,\ldots,x'_N)$
with each $x'_j$ ($1\le j\le N$) being a rational function in $c$ and $x$.
In what follows we consider $c$ to be a parameter and $e_i^c$ to be a 
rational transformation
on $\mathbb{R}^N$. We also abbreviate the symbol $\circ$ of composition. For 
instance
we write $e_i^{c_1}e_j^{c_2}$ instead of $e_i^{c_1}\circ e_j^{c_2}$.

\begin{definition} \label{def:pre-cry}
A {\em geometric pre-crystal} is a family $\B=\{x,\veps_i,\vphi_i,e_i^c\}$ 
satisfying
\begin{itemize}
\item[(i)] $e_i^{c_1}e_i^{c_2}(x)=e_i^{c_1c_2}(x),e_i^1(x)=x$,
\item[(ii)] 
$\veps_i(e_i^c(x))=c^{-1}\veps_i(x),\vphi_i(e_i^c(x))=c\vphi_i(x)$,
\item[(iii)] 
$\gamma_i(e_j^c(x))=c^{\langle\alpha^\vee_i,\alpha_j\rangle}\gamma_i(x)$.
\end{itemize}
\end{definition}

\begin{definition} \label{def:geom-cry}
A geometric pre-crystal $\B$ is called {\em geometric crystal}, if it 
further satisfies
the following relations:
\begin{equation} \label{eq:geom-cry1}
e_i^{c_1}e_j^{c_2}=e_j^{c_2}e_i^{c_1}
\end{equation}
if $\langle\alpha^\vee_i,\alpha_j\rangle=0$,
\begin{equation} \label{eq:geom-cry2}
e_i^{c_1}e_j^{c_1c_2}e_i^{c_2}=e_j^{c_2}e_i^{c_1c_2}e_j^{c_1}
\end{equation}
if 
$\langle\alpha^\vee_i,\alpha_j\rangle=\langle\alpha^\vee_j,\alpha_i\rangle=-  
1$.
\end{definition}

\begin{remark}
The definition of the geometric pre-crystal is slightly different from the 
original one
in \cite{BK}. We adopt it in order to make the formulas parallel with those 
in crystals.
\end{remark}

\begin{remark}
As mentioned in \cite{BK} Remarks after Lemma 2.1, 
\eqref{eq:geom-cry1}--\eqref{eq:geom-cry2}
can be thought of as multiplicative analogues of 
the Verma relations in the 
universal enveloping algebra $U(\geh)$. 
See also \cite{L1} Proposition 39.3.7 for the version
of the quantum enveloping algebra $U_q(\geh)$.
The relations in the cases of $\langle \alpha^\vee_i,\alpha_j \rangle
\langle \alpha^\vee_j,\alpha_i \rangle\ge2$ read as follows:
\[
e_i^{c_1}e_j^{c_1^2c_2}e_i^{c_1c_2}e_j^{c_2}=e_j^{c_2}e_i^{c_1c_2}e_j^{c_1^2  
c_2}e_i^{c_1}
\]
if 
$\langle\alpha^\vee_i,\alpha_j\rangle=-1,\langle\alpha^\vee_j,
\alpha_i\rangle=-2$,
\[
e_i^{c_1}e_j^{c_1^3c_2}e_i^{c_1^2c_2}e_j^{c_1^3c_2^2}e_i^{c_1c_2}e_j^{c_2}
=e_j^{c_2}e_i^{c_1c_2}e_j^{c_1^3c_2^2}e_i^{c_1^2c_2}
e_j^{c_1^3c_2}e_i^{c_1}  
\]
if 
$\langle\alpha^\vee_i,\alpha_j\rangle=-1,\langle\alpha^\vee_j,
\alpha_i\rangle=-3$.
\end{remark}

\begin{prop}[\cite{BK}]
Let $\B$ be a geometric crystal.
Set $s_i(x)=e_i^{\gamma_i(x)^{-1}}(x)$. Then $s_i$ ($i\in I$) generates
a birational action of the Weyl group corresponding to $\geh$.
\end{prop}

\begin{proof}
By definition, $s_i^2(x)=e_i^{c'}e_i^c(x)$ with $c=\gamma_i(x)^{-1},
c'=\gamma_i(e_i^{\gamma_i(x)^{-1}}(x))^{-1}$. From Definition 
\ref{def:pre-cry} (iii),
we have
\[
c'=((\gamma_i(x)^{-1})^2\gamma_i(x))^{-1}=\gamma_i(x).
\]
Thus $s_i^2(x)=x$ by Definition \ref{def:pre-cry} (i).

Similarly, if 
$\langle\alpha^\vee_i,\alpha_j\rangle=\langle\alpha^\vee_j,\alpha_i
\rangle=-1$,
we have
\begin{eqnarray*}
s_is_js_i(x)=e_i^{c''}e_j^{c'}e_i^c(x),\\
s_js_is_j(x)=e_j^{d''}e_i^{d'}e_j^d(x),
\end{eqnarray*}
with 
$c=d''=\gamma_i(x)^{-1},c'=d'=\gamma_i(x)^{-1}\gamma_j(x)^{-1},
c''=d=\gamma_j(x)^{-1}$.
Thus $s_is_js_i(x)=s_js_is_j(x)$ by Definition \ref{def:geom-cry}.

If $\langle\alpha^\vee_i,\alpha_j\rangle=0$, one clearly has 
$s_is_j(x)=s_js_i(x)$.
\end{proof}

\begin{remark}
This Weyl group action is an analogue of the one on crystals 
given in \cite{K2}.
\end{remark}

\begin{example} \label{ex:A}
Let $\geh=A^{(1)}_n,N=n+1,x=(x_1,x_2,\ldots,x_{n+1})$.
For $i\in I=\{0,1,\ldots,n\}$ we define
\begin{eqnarray*}
\veps_i(x)&=&x_{i+1},\quad \vphi_i(x)=x_i,\\
e^c_i(x)&=&(\ldots,cx_i,c^{-1}x_{i+1},\ldots).
\end{eqnarray*}
Here $x_0$ should be understood as $x_{n+1}$. One can check that
$\B=\{x,\veps_i,\vphi_i,e_i^c\}$ is a geometric crystal. $s_i(x)$ is given 
by
\[
s_i(x)=(\ldots,\stackrel{\scriptstyle i}{x_{i+1}},\stackrel{\scriptstyle 
i+1}{x_i},\ldots).
\]
\end{example}

\subsection{Geometric $D^{(1)}_n$-crystal $\B$} \label{subsec:defofB}

We are interested in a particular geometric crystal $\B$ of type $D^{(1)}_n$ 
given below.
Let $N=2n-1$ and consider a set of variables 
$x=(x_1,x_2,\ldots,x_n,\xb_{n-1},\ldots,\xb_1)$.
$\veps_i(x),\vphi_i(x),e_i^c(x)$ ($i\in I=\{0,1,\ldots,n\}$) are given as 
follows.

\begin{eqnarray*}
\veps_0(x)&=&x_1(\frac{x_2}{\xb_2}+1),\quad
\vphi_0(x)=\xb_1(\frac{\xb_2}{x_2}+1),\\
\veps_i(x)&=&\xb_i(\frac{x_{i+1}}{\xb_{i+1}}+1),\quad
\vphi_i(x)=x_i(\frac{\xb_{i+1}}{x_{i+1}}+1)\quad(i=1,\ldots,n-2),\\
\veps_{n-1}(x)&=&x_n\xb_{n-1},\quad \vphi_{n-1}(x)=x_{n-1},\\
\veps_n(x)&=&\xb_{n-1},\quad \vphi_n(x)=x_{n-1}x_n,\\
e^c_0(x)&=&(\xi_2^{-1}x_1,c^{-1}\xi_2x_2,\ldots,\xi_2\xb_2,c\xi_2^{-1}\xb_1)  
,\\
e^c_i(x)&=&(\ldots,c\xi_{i+1}^{-1}x_i,c^{-1}\xi_{i+1}x_{i+1},\ldots,
\xi_{i+1}\xb_{i+1},\xi_{i+1}^{-1}\xb_i,\ldots)\quad(i=1,\ldots,n-2),\\
e^c_{n-1}(x)&=&(\ldots,cx_{n-1},c^{-1}x_n,\ldots),\\
e^c_n(x)&=&(\ldots,cx_n,c^{-1}\xb_{n-1},\ldots),\\
\mbox{where }&&\xi_i=\frac{x_i+c\xb_i}{x_i+\xb_i}\quad(i=1,\ldots,n-1).
\end{eqnarray*}
We set
\[
\ell(x)=x_1x_2\cdots x_n\xb_{n-1}\xb_{n-2}\cdots\xb_1
\]
and call it the {\em level}. Note that $\ell(x)$ is invariant 
under
the transformation $e_i^c$ for any $i$ and $c$.

\begin{remark}
This $\B$ is obtained as a multiplicative analogue of the
coherent family of perfect crystals $\{B_l\}_{l\ge1}$ of type $D^{(1)}_n$ 
given in
\cite{KKM}. Of course, $e_i^c(x)$ corresponds to $\et{i}^cb$.
\end{remark}

The functions $\veps_i(x),\vphi_i(x)$ of our $\B$ satisfy the following good 
properties,
which will be used in the next section.

\begin{lemma} \label{lem:eps-phi-decomp}
\begin{itemize}

\item[(a)] If $\langle\alpha^\vee_i,\alpha_j\rangle=0$, then
\[
\veps_i(e_j^c(x))=\veps_i(x),\quad \vphi_i(e_j^c(x))=\vphi_i(x).
\]

\item[(b)] If 
$\langle\alpha^\vee_i,\alpha_j\rangle=\langle\alpha^\vee_j,
\alpha_i\rangle=-1$, then
there exist rational functions $\veps_{ij}(x),\veps_{ji}(x)$ such that
\begin{eqnarray*}
&&\veps_i(x)\veps_j(x)=\veps_{ij}(x)+\veps_{ji}(x),\\
&&\veps_{ij}(e_i^c(x))=c^{-1}\veps_{ij}(x),\quad 
\veps_{ij}(e_j^c(x))=\veps_{ij}(x),\\
&&\veps_{ji}(e_i^c(x))=\veps_{ij}(x),\quad 
\veps_{ji}(e_j^c(x))=c^{-1}\veps_{ji}(x),
\end{eqnarray*}
and $\vphi_{ij}(x),\vphi_{ji}(x)$ such that
\begin{eqnarray*}
&&\vphi_i(x)\vphi_j(x)=\vphi_{ij}(x)+\vphi_{ji}(x),\\
&&\vphi_{ij}(e_i^c(x))=c\vphi_{ij}(x),\quad 
\vphi_{ij}(e_j^c(x))=\vphi_{ij}(x),\\
&&\vphi_{ji}(e_i^c(x))=\vphi_{ij}(x),\quad 
\vphi_{ji}(e_j^c(x))=c\vphi_{ji}(x).
\end{eqnarray*}

\end{itemize}
\end{lemma}

\begin{remark} \label{rem:eps_ij}
Due to Definition \ref{def:pre-cry} (iii) it suffices to check Lemma
\ref{lem:eps-phi-decomp} only for $\vphi_i(x)$ (or $\veps_i(x)$). The 
functions
$\veps_{ij}(x),\veps_{ji}(x)$ in (b) are obtained as
\[
\veps_{ij}(x)=\vphi_{ji}(x)/(\gamma_i(x)\gamma_j(x)),\quad
\veps_{ji}(x)=\vphi_{ij}(x)/(\gamma_i(x)\gamma_j(x)),
\]
once $\vphi_{ij}(x),\vphi_{ji}(x)$ are obtained.
\end{remark}

\begin{proof}[Proof of Lemma \ref{lem:eps-phi-decomp}.]
(a) is easy. For (b) it suffices to give a list of $\vphi_{ij}(x)$ and 
$\vphi_{ji}(x)$
when 
$\langle\alpha^\vee_i,\alpha_j\rangle=\langle\alpha^\vee_j,
\alpha_i\rangle=-1$.
\begin{eqnarray*}
&&\vphi_{02}(x)=\xb_1\xb_2(\frac{\xb_3}{x_3}+1),\quad
\vphi_{20}(x)=\xb_1 x_2(\frac{\xb_3}{x_3}+1),\\
&&\vphi_{i,i+1}(x)=x_i\xb_{i+1}(\frac{\xb_{i+2}}{x_{i+2}}+1),\quad
\vphi_{i+1,i}(x)=x_i x_{i+1}(\frac{\xb_{i+2}}{x_{i+2}}+1) \quad(1\le i\le 
n-3),\\
&&\vphi_{n-2,n-1}(x)=x_{n-2}\xb_{n-1},\quad 
\vphi_{n-1,n-2}(x)=x_{n-2}x_{n-1},\\
&&\vphi_{n-2,n}(x)=x_{n-2}x_n\xb_{n-1},\quad 
\vphi_{n,n-2}(x)=x_{n-2}x_{n-1}x_n.
\end{eqnarray*}
\end{proof}

\subsection{Product}\label{subsec:product}
Let $\B_1=\{x=(x_1,\ldots,x_{N_1}),\veps_i,\vphi_i,e_i^c\},
\B_2=\{y=(y_1,\ldots,y_{N_2}),$ $\veps_i,\vphi_i,e_i^c\}$ be geometric 
pre-crystals.
For brevity we use the same symbols $\veps_i,\vphi_i,e_i^c$ for different 
geometric
pre-crystals, since they are distinguished by the sets of variables.
Now consider the set of variables 
$(x,y)=(x_1,\ldots,x_{N_1},y_1,\ldots,y_{N_2})$.
We define a new structure of a geometric pre-crystal $\B_1\times\B_2$ on 
$(x,y)$ by

\begin{eqnarray}
\veps_i(x,y)&=&\veps_i(x)+\veps_i(x)\veps_i(y)/\vphi_i(x),\\
\vphi_i(x,y)&=&\vphi_i(y)+\vphi_i(x)\vphi_i(y)/\veps_i(y),\\
e_i^c(x,y)&=&(e_i^{c_1}(x),e_i^{c_2}(y)), \label{eq:ecixy}\\
\text{where }&&c_1=\frac{c\vphi_i(x)+\veps_i(y)}{\vphi_i(x)+\veps_i(y)},\;
c_2=\frac{\vphi_i(x)+\veps_i(y)}{\vphi_i(x)+c^{-1}\veps_i(y)}.
\label{eq:c1c2}
\end{eqnarray}

\begin{remark} \label{rem:tensor}
These are analogues of formulas for the tensor product of crystals \cite{K1}
given by
\begin{eqnarray}
\veps_i(x\ot y)&=&\max(\veps_i(x),\veps_i(x)+\veps_i(y)-\vphi_i(x)),\\
\vphi_i(x\ot y)&=&\max(\vphi_i(y),\vphi_i(x)+\vphi_i(y)-\veps_i(y)),\\
\et{i}^c(x\ot y)&=&\et{i}^{c_1}(x)\ot\et{i}^{c_2}(y), \label{eq:etcixy}\\
\text{where }&&c_1=\max(c+\vphi_i(x),\veps_i(y))-\max(\vphi_i(x),\veps_i(y)),\\
&&c_2=\max(\vphi_i(x),\veps_i(y))-\max(\vphi_i(x),-c+\veps_i(y)).
\end{eqnarray}
In \eqref{eq:etcixy} $c$ can be a negative integer, in which case we understand
$\et{i}^c=\ft{i}^{-c}$.
\end{remark}

\begin{lemma}
Assume that the functions $\veps_i,\vphi_i$ for both $\B_1$ and $\B_2$ 
satisfy
the properties in Lemma \ref{lem:eps-phi-decomp}.
Then those for $\B_1\times\B_2$ also satisfy the same property.
\end{lemma}

\begin{proof}
If $\langle\alpha^\vee_i,\alpha_j\rangle=0$, the statement is clear.
Suppose 
$\langle\alpha^\vee_i,\alpha_j\rangle=\langle\alpha^\vee_j,\alpha_i\rangle
=-1$. Due to Remark \ref{rem:eps_ij} it suffices to construct desired 
functions
$\vphi_{ij}(x,y),\vphi_{ji}(x,y)$. Define
\[
\vphi_{ij}(x,y)=\vphi_{ij}(y)+\vphi_i(x)\vphi_j(y)\gamma_i(y)
+\vphi_{ij}(x)\gamma_i(y)\gamma_j(y)
\]
and $\vphi_{ji}(x,y)$ by interchanging $i$ and $j$.
To illustrate we calculate $\vphi_{ij}(e_i^c(x,y))$. Let $c_1,c_2$ be
as in \eqref{eq:c1c2}. Using Definition \ref{def:pre-cry} (iii) and
the assumption,
\begin{eqnarray*}
\vphi_{ij}(e_i^c(x,y))&=&
c_2\vphi_{ij}(y)+\frac{c_1\vphi_i(x)(c_2\vphi_{ij}(y)+\vphi_{ji}(y))}{c_2^{-  
1}\veps_i(y)}
+c_1\vphi_{ij}(x)\cdot c_2^2\gamma_i(y)\cdot c_2^{-1}\gamma_j(y)\\
&=&c_2(\vphi_{ij}(y)+\frac{c\vphi_i(x)\vphi_{ij}(y)}{\veps_i(y)})
+c(\frac{\vphi_i(x)\vphi_{ji}(y)}{\veps_i(y)}+\vphi_{ij}(x)\gamma_i(y)\gamma  
_j(y))\\
&=&c(\vphi_{ij}(y)+\frac{\vphi_i(x)\vphi_{ij}(y)}{\veps_i(y)})
+c(\frac{\vphi_i(x)\vphi_{ji}(y)}{\veps_i(y)}+\vphi_{ij}(x)\gamma_i(y)\gamma  
_j(y))\\
&=&c\vphi_{ij}(x,y).
\end{eqnarray*}

The other cases are similar.
\end{proof}

We now consider a multiple product. Let $\B_1,\ldots,\B_L$ be geometric 
pre-crystals
such that $\B_l=\{x^l=(x^l_1,\ldots,x^l_{N_l})\}$ ($l=1,\ldots,L$). Let us 
set
$\boldsymbol{x}=(x^1,\ldots,x^L)$. Then $\veps_i,\vphi_i,e^c_i$ of the 
$L$-fold product
of geometric pre-crystals $\B_1\times\cdots\times\B_L$ is given by

\begin{align}
\label{eq:vepsforxxx}
\veps_i(\boldsymbol{x}) &=
\frac{\sum_{k=1}^L \left( \prod_{j=1}^k \veps_i(x^j) \right)
\left( \prod_{j=k}^{L-1} \vphi_i(x^j) \right)}
{\prod_{j=1}^{L-1} \vphi_i(x^j)}, \\
\label{eq:vphiforxxx}
\vphi_i(\boldsymbol{x}) &=
\frac{\sum_{k=1}^L \left( \prod_{j=2}^k \veps_i(x^j) \right)
\left( \prod_{j=k}^{L} \vphi_i(x^j) \right)}
{\prod_{j=2}^{L} \veps_i(x^j)},\\
e^c_i(\boldsymbol{x}) &= (e^{c_1}_i(x^1),\ldots, e^{c_L}_i(x^L)),\\
\label{eq:def-of-c_l}
\mbox{with}\quad c_l &= \frac{\sum_{k=1}^L c^{\theta(k \leq l)}
\left( \prod_{j=2}^k \veps_i(x^j) \right)
\left( \prod_{j=k}^{L-1} \vphi_i(x^j) \right)}
{\sum_{k=1}^L c^{\theta(k \leq l-1)}
\left( \prod_{j=2}^k \veps_i(x^j) \right)
\left( \prod_{j=k}^{L-1} \vphi_i(x^j) \right)}.
\end{align}
Here $\theta(s)=1$ if $s$ is true and $=0$ otherwise.
Note that the product is associative,
$e^c_i(x^1,x^2,x^3) = 
\left(e^{c_1c_2}_i(x^1,x^2),e^{c_3}_i(x^3) \right) =
\left(e^{c_1}_i(x^1),e^{c_2c_3}_i(x^2,x^3) \right)$.
\section{Realization by unipotent matrices}\label{sec:M}

\subsection{Matrix $M(x,z)$}
\noindent
Let us introduce $2n$ by $2n$
matrices $G_i(a)$ for $0\le i\le n$ by
\begin{displaymath}
G_i(a) = E +
\begin{cases}
a (E_{i,i+1} + E_{2n-i,2n+1-i} ) & \mbox{for} \; 1 \leq i \leq n-1, \\
a (E_{n-1,n+1} + E_{n,n+2} ) & \mbox{for} \; i=n, \\
a (E_{2n-1,1} + E_{2n,2}) & \mbox{for} \; i=0.
\end{cases}
\end{displaymath}
Here $E$ is the identity matrix, $E_{ij}$ is the $(i,j)$ matrix unit,
and $a$ is a parameter.
$G_i(a)$'s have the following properties.
\begin{eqnarray}
&&G_i(a)G_i(b)=G_i(a+b),\\
&&G_i(a)G_j(b)=G_j(b)G_i(a)\quad \mbox{ if 
}\langle\alpha^\vee_i,\alpha_j\rangle=0,\\
&&G_i(a)G_j(b)G_i(c)=G_j(a')G_i(b')G_j(c')\quad \mbox{ if }
\langle\alpha^\vee_i,\alpha_j\rangle=\langle\alpha^\vee_j,\alpha_i\rangle=-1  
,
\label{eq:GiGjGi=GjGiGj}
\end{eqnarray}
where $a'=bc/(a+c),b'=a+c,c'=ab/(a+c)$.

Let $\B=\{x=(x_1,\ldots,\xb_1),\veps_i,\vphi_i,e_i^c\}$
be the geometric $D^{(1)}_n$-crystal given in the previous section.
In this subsection
we introduce a $2n$ by $2n$ matrix
$M(x,z)$ that satisfies
\begin{displaymath}
G_i \left( \frac{c-1}{z^{\delta_{i0}}\veps_i(x)} \right) M(x,z)
G_i \left( \frac{c^{-1}-1}{z^{\delta_{i0}}\vphi_i(x)} \right) =
M(e^c_i (x),z)
\end{displaymath}
for any $i$ ($0 \leq i \leq n$) and an extra parameter $z$.
We call $z$ the {\em spectral parameter}.

We give an explicit form of $M(x,z)$.
Matrix elements in the first column
$M(x,z)_{i,1}$ are given as follows.
\begin{equation}
\label{eq:firstcolumnofA}
M(x,z)_{i,1} =
\begin{cases}
x_1/\xb_1 & \mbox{for}  \; i=1, \\
x_1 \cdots x_{i-1} (1+x_i/\xb_i) & \mbox{for} \;  2 \leq i \leq n-1, \\
x_1 \cdots x_{n-1}x_n & \mbox{for} \;  i=n, \\
x_1 \cdots x_{n-1} & \mbox{for} \;  i=n+1, \\
\ell(x) / (\xb_1 \cdots \xb_{2n-i}) & \mbox{for} \;  n+2 \leq i \leq 2n.
\end{cases}
\end{equation}

Given the first column as above we define
the other matrix elements of $M(x,z)$
by the relations (\ref{eq:defMno1})-(\ref{eq:defMno6}).
\begin{equation}
\label{eq:defMno1}
M(x,z)_{i,j+1} =
\frac{M(x,z)_{i,j}}{x_j} +
\begin{cases}
(\ell (x) z - 1)/\xb_j & \mbox{for} \; i=j, \\
\ell (x) z - 1 & \mbox{for} \; i=j+1, \\
0 & \mbox{for} \; i \ne j,j+1,
\end{cases}
\end{equation}
for $1 \leq j \leq n-2$.

\begin{equation}
M(x,z)_{i,n} =
\frac{M(x,z)_{i,n-1}}{x_{n-1}} +
\begin{cases}
(\ell (x) z - 1)/\xb_{n-1} & \mbox{for} \; i=n-1, \\
\ell (x) z - 1 & \mbox{for} \; i=n+1, \\
0 & \mbox{for} \; i \ne n-1,n+1.
\end{cases}
\end{equation}

\begin{equation}
M(x,z)_{i,n+1} =
\frac{M(x,z)_{i,n-1}}{x_{n-1}x_n} +
\begin{cases}
(\ell (x) z - 1)/(\xb_{n-1}x_n) & \mbox{for} \; i=n-1, \\
\ell (x) z - 1 & \mbox{for} \; i=n, \\
0 & \mbox{for} \; i \ne n-1,n.
\end{cases}
\end{equation}

\begin{equation}
M(x,z)_{i,2n} =
z \frac{\xb_1 M(x,z)_{i,1}}{x_1} +
\begin{cases}
\frac{\xb_1}{x_1} (1 - \ell (x) z) & \mbox{for} \; i=2n, \\
z(\ell (x) z - 1) & \mbox{for} \; i=1, \\
0 & \mbox{for} \; i \ne 1,2n.
\end{cases}
\end{equation}
Let us denote $1+\xb_i/x_i$ by $\overline{X}_i$.

\begin{equation}
M(x,z)_{i,2n-1} =
z \xb_1 \overline{X}_2 M(x,z)_{i,1} +
\begin{cases}
\xb_1 \overline{X}_2 (1 - \ell (x) z) & \mbox{for} \; i=2n, \\
\frac{\xb_2}{x_2}(1 - \ell (x) z) & \mbox{for} \; i=2n-1, \\
0 & \mbox{for} \; i \ne 2n-1,2n.
\end{cases}
\end{equation}

\begin{equation}
\label{eq:defMno6}
M(x,z)_{i,2n-j} =
\xb_j \frac{\overline{X}_{j+1} M(x,z)_{i,2n+1-j}}{\overline{X}_{j} } +
\begin{cases}
\xb_j \frac{\overline{X}_{j+1}}{\overline{X}_{j} }
(1 - \ell (x) z) & \mbox{for} \; i=2n+1-j, \\
\frac{\xb_{j+1}}{x_{j+1}}(1 - \ell (x) z) & \mbox{for} \; i=2n-j, \\
0 & \mbox{for} \; i \ne 2n+1-j,2n-j,
\end{cases}
\end{equation}
for $2 \leq j \leq n-2$.
For later use we also present
matrix elements in the last row.
\begin{equation}
\label{eq:lastrowofA}
M(x,z)_{2n,i} =
\begin{cases}
\ell(x) / (x_1 \cdots x_{i-1}) & \mbox{for} \;  1 \leq i \leq n,\\
\xb_1 \cdots \xb_{n-1} & \mbox{for} \;  i=n+1, \\
\xb_1 \cdots \xb_{2n-i} (1+\xb_{2n+1-i}/x_{2n+1-i})
& \mbox{for} \;  n+2 \leq i \leq 2n-1, \\
\xb_1/x_1 & \mbox{for}  \; i=2n.
\end{cases}
\end{equation}

It can be seen that each matrix element is a polynomial in $z$ of degree at 
most $2$.
We denote the coefficients of $M(x,z)$ of degree $0,1,2$ by 
$A(x),B(x),C(x)$. Hence
we have
\[
M(x,z)=A(x)+zB(x)+z^2C(x).
\]
\begin{example}[$n=4$ case]
\label{ex:abc}
\noindent
$x=(x_1,x_2,x_3,x_4,\xb_3,\xb_2,\xb_1)$.
We use the notations
$[1]= x_1, [2]= x_2, \ldots , [\bar{1}] = \xb_1$, and
$[(2)] = 1+x_2/\xb_2, [(\bar{3})]=1+\xb_3/x_3$, etc.
For instance, $[34\bar{3}]=x_3x_4\xb_3$ and 
$[1(\bar{2})]=x_1(1+\xb_2/x_2)$.
Define $8 \times 8$ matrices $A(x), B(x), C(x)$ as follows.
\begin{displaymath}
A(x)=\left(
\begin{array}{rrrrrrrr}
\frac{[1]}{[\bar{1}]}  \\
 \left[1(2)\right]  & \frac{ [2]}{ [\bar{2}]}  \\
 \left[12(3)\right]  &  \left[2(3)\right]  & \frac{ [3]}{ [\bar{3}]}  \\
 \left[1234\right]  &  \left[234\right] &  \left[34\right] &  [4]  \\
 \left[123\right]  &  \left[23\right] &  [3] & & \frac{ 1}{ [4]}   \\
 \left[1234\bar{3}\right] &  \left[234\bar{3}\right]  & 
 \left[34\bar{3}\right] &  \left[4\bar{3}\right] &  [\bar{3}] & \frac{ 
[\bar{3}]}{ [3]}   \\
 \left[1234\bar{3}\bar{2}\right] & \left[234\bar{3}\bar{2}\right] & 
 \left[34\bar{3}\bar{2}\right] &  \left[4\bar{3}\bar{2}\right] & 
 \left[\bar{3}\bar{2}\right] &  \left[(\bar{3})\bar{2}\right] & \frac{ 
[\bar{2}]}{ [2]}  \\
 \left[1234\bar{3}\bar{2}\bar{1}\right] & 
 \left[234\bar{3}\bar{2}\bar{1}\right] & 
 \left[34\bar{3}\bar{2}\bar{1}\right] &  \left[4\bar{3}\bar{2}\bar{1}\right] 
&  \left[\bar{3}\bar{2}\bar{1}\right] & 
 \left[(\bar{3})\bar{2}\bar{1}\right] &  \left[(\bar{2})\bar{1}\right] & 
\frac{ [\bar{1}]}{ [1]}
\end{array}
\right),
\end{displaymath}
%
\begin{eqnarray*}
B(x)&=&\left(
\begin{array}{rrrrrr}
0  &  \left[1234 \bar{3}\bar{2}\right]  & \left[134 \bar{3}\bar{2}\right]  & 
\left[14\bar{3}\bar{2}\right]  & \left[1 \bar{3}\bar{2}\right]  & 
\left[1(\bar{3})\bar{2}\right] \\
 & \left[1234 \bar{3}\bar{2}\bar{1}\right] &  \left[1(2)34 
\bar{3}\bar{2}\bar{1}\right] & \left[1(2)4\bar{3}\bar{2}\bar{1}\right]  & 
\left[1(2) \bar{3}\bar{2}\bar{1}\right]  & \left[1(2)(\bar{3})\bar{2}\bar{1} 
\right] \\
&  & \left[1234 \bar{3}\bar{2}\bar{1}\right] & 
\left[12(3)4\bar{3}\bar{2}\bar{1}\right]
& \left[12(3) \bar{3}\bar{2}\bar{1}\right] & 
\left[12(3)(\bar{3})\bar{2}\bar{1}\right] \\
&  &  & & \left[1234\bar{3}\bar{2}\bar{1}\right] & 
\left[1234(\bar{3})\bar{2}\bar{1}\right] \\
&  &  & \left[1234 \bar{3}\bar{2}\bar{1}\right] & & \left[123 
(\bar{3})\bar{2}\bar{1}\right] \\
 &  &  & &   &\left[1234 \bar{3}\bar{2}\bar{1}\right] \\
 &  & &  &  &  \\
 &  &  & &  &
\end{array}
\right. \\
&&\hspace{4mm}\left.
\begin{array}{rr}
\left[1(\bar{2})\right] &   \\
\left[1(2)(\bar{2})\bar{1}\right]  &  \left[(2)\bar{1}\right] \\
\left[12(3)(\bar{2})\bar{1}\right] & \left[2(3)\bar{1}\right] \\
\left[1234(\bar{2})\bar{1}\right] & \left[234\bar{1}\right] \\
\left[123 (\bar{2})\bar{1}\right] & \left[23 \bar{1}\right] \\
\left[1234 \bar{3}(\bar{2})\bar{1}\right] & \left[234 \bar{3}\bar{1}\right] 
\\
\left[1234 \bar{3}\bar{2}\bar{1}\right] & \left[234 
\bar{3}\bar{2}\bar{1}\right] \\
 & 0
\end{array}
\right),
\quad
C(x)=\left(
\begin{array}{rrrc}
0  &  \cdots  & 0  & \left[1234 \bar{3}\bar{2}\bar{1}\right] \\
0  &  \cdots  & \cdots  & 0 \\
\vdots &  &   & \vdots \\
0 & \cdots & \cdots & 0
\end{array}
\right).
\end{eqnarray*}
The matrix $M(x,z)$ for $D_4^{(1)}$ is given by
\begin{math}
M(x,z) = A(x) + z B(x) + z^2 C(x).
\end{math}
\end{example}

Before going into the main theorem (Theorem \ref{th:GMGequalM})
we present several properties of the matrix $M(x,z) = A(x)+z B(x) + 
z^2C(x)$.
The only nonzero element of $C(x)$ is $C(x)_{1,2n}$
that is equal to $\ell(x)$.
The matrix $A(x)$ is lower triangular. 
Explicit expressions for the nonzero elements of $A(x)$ are given as follows:
\begin{itemize}
	\item $A(x)_{i,i}=x_i/\xb_i \, (1 \leq i \leq n-1), \, = x_n \, (i=n),
	\, = 1/x_n \, (i=n+1), \, =\xb_{2n+1-i}/x_{2n+1-i} \, (n+2 \leq i \leq 
	2n)$. \\
	\item For $2 \leq i \leq n-1, 1 \leq j \leq i-1$, we have
	\begin{math}
	A(x)_{i,j} = x_j \cdots x_{i-1} (1 + x_i/\xb_i )
	\end{math}
	and 
	\begin{math}
	A(x)_{2n+1-j,2n+1-i} = \xb_j \cdots \xb_{i-1} (1 + \xb_i/x_i ).
	\end{math} \\
	\item For $1 \leq j \leq n-1$ we have
	\begin{math}
	A(x)_{n,j} = x_j \cdots x_n, A(x)_{n+1,j} = x_j \cdots x_{n-1}, 
        \end{math}
        \begin{math}
	A(x)_{2n+1-j,n} = \xb_j \cdots \xb_{n-1} x_n, 
	\end{math}
	and
	\begin{math}
	A(x)_{2n+1-j,n+1} = \xb_j \cdots \xb_{n-1}.
	\end{math} \\
	\item For $1 \leq i,j \leq n-1$ we have
	\begin{math}
	A(x)_{2n+1-i,j} = (x_j \cdots x_n ) \times (\xb_i \cdots \xb_{n-1}).
	\end{math}
\end{itemize}
\begin{lemma}
Let $F_i(a)={}^t\!G_i(a)$ and
\begin{displaymath}
d(x) = \mbox{diag} \left( x_1/\xb_1, \ldots, x_{n-1}/\xb_{n-1},x_n,
1/x_n,\xb_{n-1}/x_{n-1},\ldots,\xb_1/x_1 \right).
\end{displaymath}
Then the matrix $A(x)$ has the following factorization.
\begin{eqnarray}
\label{eq:factorizationofA}
A(x) &=& F_1(\xb_1) F_2(\xb_2) \cdots F_{n-2}(\xb_{n-2}) F_n(\xb_{n-1})
d(x) \\
& \times & F_{n-1}(x_{n-1}) F_{n-2}(x_{n-2}) \cdots F_2(x_2)F_1(x_1)
\nonumber
\end{eqnarray}
\end{lemma}
\begin{proof}
$F_1(\xb_1)^{-1} A(x)F_1(x_1)^{-1}$ is given by deleting the elements in
the 1st column and the $2n$-th row of $A(x)$ except the diagonal elements.
To illustrate we take the $n=4$ case. See Example \ref{ex:abc}.
\begin{displaymath}
F_1(\xb_1)^{-1} A(x)F_1(x_1)^{-1}=\left(
\begin{array}{rrrrrrrr}
\frac{[1]}{[\bar{1}]}  \\
  & \frac{ [2]}{ [\bar{2}]}  \\
  &  \left[2(3)\right]  & \frac{ [3]}{ [\bar{3}]}  \\
  &  \left[234\right] &  \left[34\right] &  [4]  \\
  &  \left[23\right] &  [3] & & \frac{ 1}{ [4]}   \\
 &  \left[234\bar{3}\right]  &  \left[34\bar{3}\right] & 
 \left[4\bar{3}\right] &  [\bar{3}] & \frac{ [\bar{3}]}{ [3]}   \\
 & \left[234\bar{3}\bar{2}\right] &  \left[34\bar{3}\bar{2}\right] & 
 \left[4\bar{3}\bar{2}\right] &  \left[\bar{3}\bar{2}\right] & 
 \left[(\bar{3})\bar{2}\right] & \frac{ [\bar{2}]}{ [2]}  \\
 &   &   &   &   &   &   & \frac{ [\bar{1}]}{ [1]}
\end{array}
\right).
\end{displaymath}
Similarly, $F_2(\xb_2)^{-1} F_1(\xb_1)^{-1} A(x)
F_1(x_1)^{-1} F_2(x_2)^{-1}$ is given by deleting the elements in
the 2nd column and the $(2n-1)$-th row of
$F_1(\xb_1)^{-1} A(x)F_1(x_1)^{-1}$ except the diagonal elements.
This process continues until we obtain $d(x)$.
For example, the final step for the $n=4$ case goes as
\begin{align*}
\left(
\begin{array}{cccc}
1 & & & \\
  & 1 & & \\
-\left[\bar{3}\right] & & 1 & \\
 & -\left[\bar{3}\right] & & 1
\end{array}
\right)
&\left(
\begin{array}{rrrr}
\frac{[3]}{ [\bar{3}]}  \\
\left[34\right] &  [4]  \\
\left[3\right] & & \frac{ 1}{ [4]}   \\
\left[34\bar{3}\right] &  \left[4\bar{3}\right] &  [\bar{3}] &
\frac{ [\bar{3}]}{ [3]}
\end{array}
\right)
\left(
\begin{array}{cccc}
1 & & & \\
-\left[3\right]   & 1 & & \\
& & 1 & \\
& & -\left[3\right] & 1
\end{array}
\right) \\
& =
\left(
\begin{array}{rrrr}
\frac{[3]}{ [\bar{3}]}  \\
&  [4]  \\
& & \frac{ 1}{ [4]}   \\
& & & \frac{ [\bar{3}]}{ [3]}
\end{array}
\right).
\end{align*}
\end{proof}
\begin{remark}
We note that the sequence $(1,2,\ldots,n-2,n,n-1,\ldots,1)$ 
extracted from the indices of $F$ in the right hand side of 
\eqref{eq:factorizationofA} is essentially the same as the 
one given in \cite{KMOTU,HKT2}.
\end{remark}
\begin{prop} \label{prop:meqdpd}
\begin{equation}
\label{eq:meqdpd}
M(x,\ell(x)^{-1}) = \frac{1}{\ell (x)} {\mathcal D}_2 (x) P
{\mathcal D}_1 (x),
\end{equation}
where
\begin{eqnarray*}
{\mathcal D}_1 (x) &=& \mbox{diag} \left(
A(x)_{2n,1},A(x)_{2n,2},\ldots,A(x)_{2n,2n} \right),\\
{\mathcal D}_2 (x) &=& \mbox{diag} \left(
A(x)_{1,1},A(x)_{2,1},\ldots,A(x)_{2n,1} \right),
\end{eqnarray*}
and $P$ is the matrix with all the entries being 1.
\end{prop}

\begin{corollary} \label{cor:B}
\[
B(x)={\mathcal D}_2(x)P{\mathcal D}_1(x)-\ell(x)A(x)-\ell(x)^{-1}C(x).
\]
\end{corollary}

Set $T=(\delta_{i+j,n+1}(-1)^{i-1})_{1 \leq i,j \leq n}$
and define a $2n$ by $2n$ matrix $S$ by
\begin{equation}
\label{eq:signmatrix}
S = \left(
\begin{array}{cc}
   & T \\
{}^t T &
\end{array}
\right).
\end{equation}
It is clear that $S^{-1} = \, ^t \! S = S$.
\begin{prop}\label{pr:MSMS}
$M(x,z) S \, ^t \! M(x,z) S = (1-z \ell(x))^2 E$.
\end{prop}
\begin{proof}
For brevity we write $M(x,z)=M$, $A(x)=A$, $B(x)=B$,
$C(x)=C$, ${\mathcal D}_1 (x)={\mathcal D}_1$,
${\mathcal D}_2 (x) = {\mathcal D}_2$, and $\ell(x)=\ell$.
For any $2n$ by $2n$ matrix $X$ we denote $S \, ^t \! X S$
by $\check{X}$.
Then we have
\begin{displaymath}
M \check{M} = A \check{A} + z (A \check{B} + B \check{A}) +
z^2 (A \check{C} + B \check{B} + C \check{A}).
\end{displaymath}
First we see that $A \check{A} = E$ since $A$ has the factorization
(\ref{eq:factorizationofA}) and any factor $X$ thereof satisfies
$X \check{X} = E$.
By Corollary \ref{cor:B} we have
\begin{math}
B = {\mathcal D}_2 P {\mathcal D}_1 - \ell A - \ell^{-1} C.
\end{math}
Let us denote
$\ell^{-1} C =(\delta_{i,1} \delta_{j,2n})_{1 \leq i,j \leq 2n}$ by $H$.
Then we have
\begin{align*}
A \check{B} &= A \check{{\mathcal D}_1} \check{P} \check{{\mathcal D}_2}
- \ell E - A H, \\
B \check{A} &= {\mathcal D}_2 P {\mathcal D}_1 \check{A}
- \ell E - H \check{A}.
\end{align*}
Here we have used $A \check{A} = E$ and $\check{H}=H$.
It is easy to see that
$A \check{{\mathcal D}_1} \check{P} \check{{\mathcal D}_2}=H \check{A} $,
${\mathcal D}_2 P {\mathcal D}_1 \check{A} = A H$,
and thereby $A \check{B} + B \check{A} = -2 \ell E$.

It remains to check the desired relation for one special value
of $z$ other than $z=0$.
Let us put $z=\ell^{-1}$.
By Proposition \ref{prop:meqdpd} we have
\begin{math}
M \check{M} = \ell^{-2} {\mathcal D}_2 P {\mathcal D}_1
\check{{\mathcal D}_1} \check{P} \check{{\mathcal D}_2}.
\end{math}
By using $A \check{A}=E$ we see
\begin{math}
P {\mathcal D}_1 \check{{\mathcal D}_1} \check{P} = O,
\end{math}
where $O$ is the zero matrix.
The proof is completed.
\end{proof}
\begin{corollary}\label{coro:detM}
$\det M(x,z) = (1-z \ell(x))^{2n}$.
\end{corollary}

On the set of variables of the geometric crystal $\B$ we introduce involutive 
automorphisms $\sigma_1$, $\sigma_n$ and $\tau$ by
\begin{align}
\label{eq:sigma1o}
\sigma_1 &: x_1 \longleftrightarrow \xb_1, \\
\label{eq:sigmano}
\sigma_n &: x_{n-1} \rightarrow x_{n-1}x_n, \quad
             \xb_{n-1} \rightarrow \xb_{n-1}x_n, \quad
             x_n \rightarrow 1/x_n,\\
\label{eq:tauo}
\tau &: x_i \longleftrightarrow \xb_i, \quad
        \xb_i \longleftrightarrow x_i, \quad (1 \leq i \leq n-2) \\
     &\;\;\; x_{n-1} \rightarrow \xb_{n-1}x_n, \quad
             \xb_{n-1} \rightarrow x_{n-1}x_n, \quad
             x_n \rightarrow 1/x_n. \nonumber
\end{align}
Note that $\sigma_1$, $\sigma_n$ and $\tau$ are mutually commutative. 
For simplicity we write $x^{\sigma_1}$ to mean $\sigma_1(x)$, etc. 
We can easily check

\begin{lemma} \label{lem:sigma}
\begin{align}
\label{eq:sigma_1}
&\veps_1(x^{\sigma_1})=\veps_0(x),\quad \vphi_1(x^{\sigma_1})=
\vphi_0(x),\quad e_1^c(x^{\sigma_1})=(e^c_0(x))^{\sigma_1},\\
\label{eq:sigma_n}
&\veps_{n-1}(x^{\sigma_n})=\veps_n(x),\quad 
\vphi_{n-1}(x^{\sigma_n})=\vphi_n(x),\quad
e_{n-1}^c(x^{\sigma_n})=(e^c_n(x))^{\sigma_n}.
\end{align}
The other $\veps_i$ and $\vphi_i$ are invariant under $\sigma_1$ and 
$\sigma_n$.
\end{lemma}

Let us introduce the $2n$ by $2n$ matrices:
\begin{align}
\label{eq:JMatrixA}
J_1(z) &= z E_{1,2n} + z^{-1} E_{2n,1} + {\textstyle \sum_{i=2}^{2n-1}}
 E_{i,i}, \\
 \label{eq:JMatrixB}
J_n &= E_{n,n+1} + E_{n+1,n} + {\textstyle \sum_{i=1}^{n-1}}
(E_{i,i}+E_{2n+1-i,2n+1-i}), \\
\label{eq:JMatrixC}
J &= {\textstyle \sum_{i=1}^{2n}} E_{i,2n+1-i}.
\end{align}
The following lemma is immediate.

\begin{lemma} \label{lem:J}
\begin{align}
\label{eq:J_1}
J_1(z)G_1(a)&=G_0(\frac{a}{z})J_1(z),\\
\label{eq:J_n}
J_nG_{n-1}(a)&=G_n(a)J_n.
\end{align}
\end{lemma}
\begin{lemma}
\label{lem:JMJequalsM}
\begin{align}
\label{eq:JMJequalsMno1}
J_1(z) M(x,z) J_1(z) &= M(x^{\sigma_1},z), \\
\label{eq:JMJequalsMno2}
J_n M(x,z) J_n &= M(x^{\sigma_n},z), \\
\label{eq:JMJequalsMno3}
J \,^t \!M(x,z) J &= M(\tau(x),z).
\end{align}
\end{lemma}
\begin{proof}
First we prove (\ref{eq:JMJequalsMno1}).
In this case, only
the first and $2n$-th rows and columns
are changed.
Let us consider $J_1(z) A(x) J_1(z)$.
All the elements in the first column and the $2n$-th row of
this matrix are zero except the $(1,1)$ and
$(2n,2n)$ elements.
And we have
\begin{align*}
\left( J_1(z) A(x) J_1(z) \right)_{1,i} &=
\begin{cases}
\xb_1/x_1 & \mbox{for}\,i=1,\\
\ell(x)z^2 & \mbox{for}\,i=2n,\\
A(x)_{2n,i}z & \mbox{for}\,i \ne 1,2n,
\end{cases} \\
\left( J_1(z) A(x) J_1(z) \right)_{i,2n} &=
\begin{cases}
\ell(x)z^2 & \mbox{for}\,i=1,\\
x_1/\xb_1 & \mbox{for}\,i=2n,\\
A(x)_{i,1}z & \mbox{for}\,i \ne 1,2n.
\end{cases} \\
\end{align*}
Note that $\xb_1/x_1 = A(x^{\sigma_1})_{1,1}$,
$A(x)_{2n,i} = B(x^{\sigma_1})_{1,i} \, (i \ne 1,2n)$,
$x_1/\xb_1 = A(x^{\sigma_1})_{2n,2n} $, and
$A(x)_{i,1} = B(x^{\sigma_1})_{i,2n} \, (i \ne 1,2n)$.
Here we have used
\begin{align*}
B(x)_{1,i} &=
\begin{cases}
0 & \mbox{for}\,i=1,2n,\\
A(x^{\sigma_1})_{2n,i} & \mbox{for}\,i \ne 1,2n,
\end{cases} \\
B(x)_{i,2n} &=
\begin{cases}
0 & \mbox{for}\,i=1,2n,\\
A(x^{\sigma_1})_{i,1} & \mbox{for}\,i \ne 1,2n,
\end{cases} \\
\end{align*}
which is obtained from Corollary \ref{cor:B}.
Let us consider $J_1(z) B(x) J_1(z)$.
{}From Corollary \ref{cor:B} we find
$B(x)_{2n,i} = B(x)_{i,1} =0$ for $1 \leq i \leq 2n$.
Therefore all elements in the first row and the $2n$-th column of
this matrix are zero.
And we have
\begin{align*}
\left( J_1(z) B(x) J_1(z) \right)_{2n,i} &=
\begin{cases}
0 & \mbox{for}\,i=1,2n,\\
B(x)_{1,i}z^{-1} & \mbox{for}\,i \ne 1,2n,
\end{cases} \\
\left( J_1(z) B(x) J_1(z) \right)_{i,1} &=
\begin{cases}
0 & \mbox{for}\,i=1,2n,\\
B(x)_{i,2n} z^{-1} & \mbox{for}\,i \ne 1,2n.
\end{cases} \\
\end{align*}
Note that $B(x)_{1,i}=A(x^{\sigma_1})_{2n,i}$ and
$B(x)_{i,2n} = A(x^{\sigma_1})_{i,1}$ for $i \ne 1,2n$.
Let us consider $J_1(z) C(x) J_1(z)$.
The only nonzero element of this matrix is the $(2n,1)$ element
that takes the value $\ell(x) z^{-2}$.
Putting all the above relations together we obtain 
(\ref{eq:JMJequalsMno1}).

Next let us consider (\ref{eq:JMJequalsMno2}).
{}From (\ref{eq:factorizationofA}) we have
\begin{eqnarray*}
J_n A(x)J_n &=& F_1(\xb_1) F_2(\xb_2) \cdots F_{n-2}(\xb_{n-2})
F_{n-1}(\xb_{n-1})
d(x^{\sigma_n}) \\
& \times & F_{n}(x_{n-1}) F_{n-2}(x_{n-2}) \cdots F_2(x_2)F_1(x_1).
\end{eqnarray*}
Since $F_{n-1}(\xb_{n-1}) d(x^{\sigma_n})=
d(x^{\sigma_n}) F_{n-1}(x_{n-1}x_n)$
and $d(x^{\sigma_n}) F_{n}(x_{n-1}) = F_{n}($ 
\par\noindent
$\xb_{n-1}x_n) d(x^{\sigma_n})$ 
we obtain $J_n A(x)J_n = A(x^{\sigma_n})$.
Note that $F_{n-1}$ and $F_n$ are commutative.
Clearly we have
$J_n C(x)J_n = C(x) = C(x^{\sigma_n})$.
Then by Corollary \ref{cor:B} we obtain
$J_n B(x)J_n = B(x^{\sigma_n})$.
The relation (\ref{eq:JMJequalsMno2}) is proved.

Finally we show (\ref{eq:JMJequalsMno3}).
{}From (\ref{eq:factorizationofA}) we have
\begin{eqnarray*}
J  \,^t \!A(x)J &=& F_1(x_1) F_2(x_2) \cdots F_{n-2}(x_{n-2})
F_{n-1}(x_{n-1})
d(\tau(x)) \\
& \times & F_{n}(\xb_{n-1}) F_{n-2}(\xb_{n-2}) \cdots 
F_2(\xb_2)F_1(\xb_1).
\end{eqnarray*}
Since $F_{n-1}(x_{n-1}) d(\tau(x))=
d(\tau(x)) F_{n-1}(\xb_{n-1}x_n)$
and $d(\tau(x)) F_{n}(\xb_{n-1}) = F_{n}(x_{n-1}$
\par\noindent
$x_n) d(\tau(x))$
we obtain
$J \,^t \!A(x)J = A(\tau(x))$.
Clearly we have
$J \,^t \!C(x)J = C(x) = C(\tau(x))$.
Then by Corollary \ref{cor:B} we can obtain
$J \,^t \!B(x)J = B(\tau(x))$.
The relation (\ref{eq:JMJequalsMno3}) is proved.
\end{proof}

The following theorem states that the transformation $e^c_i$ of the 
geometric
$D^{(1)}_n$-crystal $\B$ is realized as the multiplications of unipotent 
matrices.
\begin{theorem}
\label{th:GMGequalM}
For $0 \leq i \leq n$ we have
\begin{equation}
G_i \left( \frac{c-1}{z^{\delta_{i0}}\veps_i(x)} \right) M(x,z)
G_i \left( \frac{c^{-1}-1}{z^{\delta_{i0}}\vphi_i(x)} \right) =
M(e^c_i (x),z).
\end{equation}
\end{theorem}

\begin{proof}
First note that the $i=n$ case can be derived from the $i=n-1$ case as 
follows.
\begin{align*}
M(e^c_n (x),z) &= J_n M(e^c_{n-1} (x^{\sigma_n}),z) J_n \\
&= J_n G_{n-1} \left( \frac{c-1}{\veps_{n-1}(x^{\sigma_n})}
\right) M(x^{\sigma_n},z)\;
G_{n-1} \left( \frac{c^{-1}-1}{\vphi_{n-1}(x^{\sigma_n})} \right) J_n \\
%
%
&= G_{n} \left( \frac{c-1}{\veps_{n}(x)}
\right) J_n M(x^{\sigma_n},z) J_n\;
G_{n} \left( \frac{c^{-1}-1}{\vphi_{n}(x)} \right) \\
&= G_{n} \left( \frac{c-1}{\veps_{n}(x)}
\right) M(x,z)\;
G_{n} \left( \frac{c^{-1}-1}{\vphi_{n}(x)} \right).
\end{align*}
Here we have used 
\eqref{eq:sigma_n},\eqref{eq:J_n},
\eqref{eq:JMJequalsMno2}.
In a similar way the $i=0$ case can be derived from the $i=1$ case 
by using \eqref{eq:sigma_1},\eqref{eq:J_1},\eqref{eq:JMJequalsMno1}.
Thus we can assume $1 \leq i \leq n-1$.
The identity to be proved is of degree 2 with respect to $z$.
However,
it is easy to see that the coefficients of $z^2$ in the both hand sides 
coincide.
Thus it remains to check the identity with two distinct values of $z$,
which we take to be $\ell(x)^{-1}$ and $0$.

First let us put $z = \ell(x)^{-1}$. From Proposition \ref{prop:meqdpd}
the identity to be proved reads as
\begin{equation}
G_i \left( \frac{c-1}{\veps_i(x)} \right)
{\mathcal D}_2(x) P {\mathcal D}_1(x)
G_i \left( \frac{c^{-1}-1}{\vphi_i(x)} \right) =
{\mathcal D}_2(e^c_i (x)) P {\mathcal D}_1(e^c_i (x)).
\end{equation}
This equality indeed holds since we have
\begin{eqnarray}
G_i \left( \frac{c-1}{\veps_i(x)} \right)
{\mathcal D}_2(x) P &=&
{\mathcal D}_2(e^c_i (x)) P, \\
P {\mathcal D}_1(x)
G_i \left( \frac{c^{-1}-1}{\vphi_i(x)} \right) &=&
P {\mathcal D}_1(e^c_i (x)),
\end{eqnarray}
or equivalently
\begin{eqnarray*}
&&A(e^c_i (x))_{k,1} = \begin{cases}
A(x)_{k,1} + A(x)_{k+1,1} \frac{c-1}{\veps_i(x)}
& \mbox{for}\; k=i, 2n-i,\\
A(x)_{k,1} & \mbox{otherwise},
\end{cases} \\
&&A(e^c_i (x))_{2n,k} = \begin{cases}
A(x)_{2n,k-1} \frac{c^{-1}-1}{\vphi_i(x)}+ A(x)_{2n,k}
& \mbox{for}\; k=i+1, 2n+1-i,\\
A(x)_{2n,k} & \mbox{otherwise}.
\end{cases}
\end{eqnarray*}
One can check them from (\ref{eq:firstcolumnofA}) and 
(\ref{eq:lastrowofA}).

Next let us put $z=0$.
We consider the $i=1$ case. The other cases are similar.
The identity to be proved reads as
\begin{equation}
\label{eq:GAGequalA}
G_1 \left( \frac{c-1}{\veps_1(x)} \right) A(x)
G_1 \left( \frac{c^{-1}-1}{\vphi_1(x)} \right) =
A(e^c_1 (x)).
\end{equation}
Substitute the factorized form of $A(x)$
(\ref{eq:factorizationofA}) into the left hand side.
Exchange the leftmost two matrices as
\begin{displaymath}
G_1 \left( \frac{c-1}{\veps_1(x)} \right) F_1(\xb_1) =
F_1(\xi_2^{-1} \xb_1) \tilde{G}_1 \left( \frac{c-1}{\veps_1(x)} \right).
\end{displaymath}
Here $\tilde{G}_1(a)$ is the matrix whose 1st and $(2n-1)$-th (resp.~2nd 
and
$2n$-th) diagonal elements are $\xi_2$ (resp.~$\xi_2^{-1}$) and the other
matrix elements are the same as those of $G_1(a)$.
Recall $\xi_2$ is defined in section \ref{subsec:defofB}.
Then the next step is to exchange the second and third matrices as
\begin{displaymath}
\tilde{G}_1 \left( \frac{c-1}{\veps_1(x)} \right) F_2(\xb_2) =
F_2(\xi_2 \xb_2) \tilde{G}_1 \left( \frac{c-1}{\veps_1(x)} \right).
\end{displaymath}
Similarly, for the rightmost three matrices we have
\begin{displaymath}
F_2(x_2)F_1(x_1) G_1 \left( \frac{c^{-1}-1}{\vphi_1(x)} \right) =
\hat{G}_1 \left( \frac{c^{-1}-1}{\vphi_1(x)} \right) F_2(c^{-1}\xi_2 x_2)
F_1(c \xi_2^{-1} x_1).
\end{displaymath}
Here $\hat{G}_1(a)$ is the matrix whose 1st and $(2n-1)$-th (resp.~2nd and
$2n$-th) diagonal elements are $c \xi_2^{-1}$ (resp.~$c^{-1} \xi_2$) and the 
other
matrix elements are the same as those of $G_1(a)$.
The other $F_i$'s commute with $\tilde{G}_1$ and $\hat{G}_1$.
Thus the left hand side of (\ref{eq:GAGequalA}) leads to
\begin{displaymath}
F_1(\xi_2^{-1} \xb_1) F_2(\xi_2 \xb_2) \cdots
\tilde{G}_1 \left( \frac{c-1}{\veps_1(x)} \right) d(x)
\hat{G}_1 \left( \frac{c^{-1}-1}{\vphi_1(x)} \right) \cdots
F_2(c^{-1} \xi_2 x_2) F_1(c \xi_2^{-1} x_1).
\end{displaymath}
Thanks to
\begin{displaymath}
d(x) \hat{G}_1 \left( \frac{c^{-1}-1}{\vphi_1(x)} \right) =
\tilde{G}_1 \left( \frac{c-1}{\veps_1(x)} \right)^{-1} d(e^c_1(x)),
\end{displaymath}
we obtain the right hand side of \eqref{eq:GAGequalA}.

The proof is completed.
\end{proof}

\subsection{Product of $M(x,z)$}

Throughout this subsection we consider the $L$-fold product $\B^{\times L}$ 
of
the geometric $D^{(1)}_n$-crystal $\B$ given in section 
\ref{subsec:defofB}.
Our main theorem here is Theorem \ref{th:ijieqjij} which states that
$\B^{\times L}$ also becomes a geometric crystal.

\begin{theorem}
\label{th:MMMequalGMMMG}
Let $\boldsymbol{x} = (x^1,\ldots,x^L)$ the set of variables of $\B^{\times 
L}$ and
set $e^c_i(\boldsymbol{x})=(y^1,\cdots,y^L)$.
Then we have
\begin{displaymath}
M(y^1,z) \cdots M(y^L,z) =
G_i\left(\frac{c-1}{z^{\delta_{i0}}\veps_i(\boldsymbol{x})} \right)
M(x^1,z) \cdots M(x^L,z) \;
G_i\left(\frac{c^{-1}-1}{z^{\delta_{i0}}\vphi_i(\boldsymbol{x})} \right),
\end{displaymath}
for $0 \leq i \leq n$.
\end{theorem}
\begin{proof}
{}From Theorem \ref{th:GMGequalM} we have
\begin{displaymath}
M(y^l,z) = M(e^{c_l}_i(x^l),z) =
G_i\left(\frac{c_l-1}{z^{\delta_{i0}}\veps_i(x^l)} \right)
M(x^l,z) \;
G_i\left(\frac{c_l^{-1}-1}{z^{\delta_{i0}}\vphi_i(x^l)} \right),
\end{displaymath}
for $1 \leq l \leq L$.
Using \eqref{eq:vepsforxxx}--\eqref{eq:def-of-c_l} we can check
\begin{eqnarray*}
&& \frac{c_{l-1}^{-1}-1}{\vphi_i(x^{l-1})} +
\frac{c_l-1}{\veps_i(x^l)} = 0 \quad \mbox{for $2 \leq l \leq L$}, \\
&&\frac{c_1-1}{\veps_i(x^1)} =
\frac{c-1}{\veps_i(\boldsymbol{x})}, \quad
\frac{c_L^{-1}-1}{\vphi_i(x^L)}=
\frac{c^{-1}-1}{\vphi_i(\boldsymbol{x})},
\end{eqnarray*}
which finish the proof.
\end{proof}

The following theorem plays a crucial role in the subsequent part of this 
paper.
It asserts that we can retrieve components $x^l$ from the product of 
matrices
$M(x^1,z)\cdots M(x^L,z)$.

\begin{theorem}
\label{th:uniqueness}
Let $\boldsymbol{x} = (x^1,\ldots,x^L)$ and
$\boldsymbol{y} = (y^1,\ldots,y^L)$ be sets of variables of $\B^{\times 
L}$.
Suppose
\begin{eqnarray*}
\ell(x^i)&=&\ell(y^i)\quad(1\le i\le L),\\
\ell(x^i)&\neq&\ell(x^j)\quad (i\neq j),
\end{eqnarray*}
and
\begin{equation}
\label{eq:mdotsmeqmdosm}
M(x^{1},z) \cdots M(x^{L},z) =M(y^{1},z) \cdots M(y^{L},z).
\end{equation}
Then we have $x^{i} = y^{i}\;(1 \leq i \leq L)$.
\end{theorem}
\begin{proof}
Denote $\ell(x^{1}) = \ell(y^{1})$ by $a^{-1}$.
Let us put $z=a$ in (\ref{eq:mdotsmeqmdosm}). Using Proposition 
\ref{prop:meqdpd}
we obtain
\begin{align*}
&{\mathcal D}_2 (y^{1})^{-1} {\mathcal D}_2 (x^{1}) P \\
&=P {\mathcal D}_1 (y^{1}) M(y^{2},a) \cdots M(y^{L},a)
M(x^{L},a)^{-1} \cdots M(x^{2},a)^{-1} {\mathcal D}_1 (x^{1})^{-1}.
\end{align*}
Here we used the property that $\det M(x,z) \ne 0$ unless 
$z=\ell(x)^{-1}$.
See Corollary \ref{coro:detM}.
We see that the both sides of this equation should be a scalar multiple of
$P$.
{}From the $2n$-th row of the left hand side
we see that the scalar should be $1$. Hence we have
${\mathcal D}_2 (x^{1}) = {\mathcal D}_2 (y^{1})$,
that implies $x^{1} = y^{1}$.
In the same way we can deduce $x^{i} = y^{i}$ for all $i$.
\end{proof}

We prepare lemmas. Let $\boldsymbol{x} = (x^1,\ldots,x^L)$ be the set of 
variables
of $\B^{\times L}$. Define $\boldsymbol{x}^{\sigma_a}
=((x^1)^{\sigma_a},\ldots,(x^L)^{\sigma_a})$ for $a=1,n$.

\begin{lemma}
\label{lem:vepsvphisigma}
\begin{align*}
&\veps_{1}(\boldsymbol{x}^{\sigma_1})
=\veps_0(\boldsymbol{x}), \quad
\vphi_{1}(\boldsymbol{x}^{\sigma_1})
=\vphi_0(\boldsymbol{x}),\\
&\veps_{n-1}(\boldsymbol{x}^{\sigma_n})
=\veps_n(\boldsymbol{x}), \quad
\vphi_{n-1}(\boldsymbol{x}^{\sigma_n})
=\vphi_n(\boldsymbol{x}).
\end{align*}
\end{lemma}

\begin{proof}
The $L=1$ case is given in Lemma \ref{lem:sigma}.
By (\ref{eq:vepsforxxx}) and (\ref{eq:vphiforxxx}) we see that
the claim is also true for $L \geq 2$.
\end{proof}

\begin{lemma}
\label{lem:Dynkinrightmost}
\begin{align*}
e^c_{1} (\boldsymbol{x}^{\sigma_1}) &=
(e^c_0 (\boldsymbol{x}))^{\sigma_1}, \\
e^c_{n-1} (\boldsymbol{x}^{\sigma_n}) &=
(e^c_n (\boldsymbol{x}))^{\sigma_n}.
\end{align*}
\end{lemma}

\begin{proof}
Denote the left (resp.~right) hand side of the 
second relation by
$\boldsymbol{y}=(y^1,\cdots,y^L)$
(resp.~$\boldsymbol{w}=(w^1,\cdots,w^L)$).
We have
\begin{align*}
&M(y^1,z) \cdots M(y^L,z) \\
&=
G_{n-1}\left(\frac{c-1}{\veps_{n-1}(\boldsymbol{x}^{\sigma_n})} \right)
M((x^1)^{\sigma_n},z) \cdots M((x^L)^{\sigma_n},z) \;
G_{n-1}\left(\frac{c^{-1}-1}{\vphi_{n-1}(\boldsymbol{x}^{\sigma_n})}
\right) \\
&=
G_{n-1}\left(\frac{c-1}{\veps_{n}(\boldsymbol{x})} \right)
J_n M(x^1,z) \cdots M(x^L,z) J_n\;
G_{n-1}\left(\frac{c^{-1}-1}{\vphi_{n}(\boldsymbol{x})}
\right) \\
&=
J_n \; G_{n}\left(\frac{c-1}{\veps_{n}(\boldsymbol{x})} \right)
M(x^1,z) \cdots M(x^L,z) \;
G_{n}\left(\frac{c^{-1}-1}{\vphi_{n}(\boldsymbol{x})}
\right) J_n\\
&= M(w^1,z) \cdots M(w^L,z) .
\end{align*}
Here we have used Lemma \ref{lem:J}, \ref{lem:JMJequalsM}, 
\ref{lem:vepsvphisigma}
and Theorem \ref{th:MMMequalGMMMG}.
Now Theorem \ref{th:uniqueness} tells 
$\boldsymbol{y}=\boldsymbol{w}$.
The first relation is similar.
\end{proof}

Our main theorem in this section is

\begin{theorem}
\label{th:ijieqjij}
$\B^{\times L}$ is a geometric crystal. Namely,
\begin{equation}
\label{eq:ijeqji}
e^d_j e^c_i (\boldsymbol{x}) =
e^c_i e^d_j (\boldsymbol{x})
\end{equation}
for $\langle\alpha^\vee_i,\alpha_j\rangle=0$, and
\begin{equation}
\label{eq:ijieqjij}
e^d_i e^{cd}_j e^c_i (\boldsymbol{x}) =
e^c_j e^{cd}_i e^d_j (\boldsymbol{x})
\end{equation}
for 
$\langle\alpha^\vee_i,\alpha_j\rangle=\langle\alpha^\vee_j,\alpha_i\rangle=-  
1$.
\end{theorem}
\begin{proof}
We prove (\ref{eq:ijieqjij}) first.
Denote the left (resp.~right) hand side of (\ref{eq:ijieqjij}) by
$\boldsymbol{y}=(y^1,\cdots,y^L)$
(resp.~$\boldsymbol{w}=(w^1,\cdots,w^L)$).
It suffices to prove
$M(y^1,z) \cdots M(y^L,z) = M(w^1,z) \cdots M(w^L,z)$
by Theorem \ref{th:uniqueness}.
Let us define
\begin{align}
\label{eq:PPQQ1}
P &= G_i\left(\frac{d-1}{\veps_{i}(e^{cd}_j e^c_i(
\boldsymbol{x}))} \right)
G_j\left(\frac{cd-1}{\veps_{j}(e^c_i(
\boldsymbol{x}))} \right)
G_i\left(\frac{c-1}{\veps_{i}(\boldsymbol{x})} \right), \\
P' &= G_j\left(\frac{c-1}{\veps_{j}(e^{cd}_i e^d_j(
\boldsymbol{x}))} \right)
G_i\left(\frac{cd-1}{\veps_{i}(e^d_j(
\boldsymbol{x}))} \right)
G_j\left(\frac{d-1}{\veps_{j}(\boldsymbol{x})} \right), \\
Q &= G_i\left(-\frac{d^{-1}-1}{\vphi_{i}(e^{cd}_j e^c_i(
\boldsymbol{x}))} \right)
G_j\left(-\frac{(cd)^{-1}-1}{\vphi_{j}(e^c_i(
\boldsymbol{x}))} \right)
G_i\left(-\frac{c^{-1}-1}{\vphi_{i}(\boldsymbol{x})} \right), \\
\label{eq:PPQQ4}
Q' &= G_j\left(-\frac{c^{-1}-1}{\vphi_{j}(e^{cd}_i e^d_j(
\boldsymbol{x}))} \right)
G_i\left(-\frac{(cd)^{-1}-1}{\vphi_{i}(e^d_j(
\boldsymbol{x}))} \right)
G_j\left(-\frac{d^{-1}-1}{\vphi_{j}(\boldsymbol{x})} \right).
\end{align}
For notational simplicity we have supposed $i,j \ne 0$.
Otherwise the associated arguments should be divided by $z$.
{}From Theorem \ref{th:MMMequalGMMMG} we have
\begin{align*}
M(y^1,z) \cdots M(y^L,z) &= P M(x^1,z) \cdots M(x^L,z) Q^{-1},\\
M(w^1,z) \cdots M(w^L,z) &= P' M(x^1,z) \cdots M(x^L,z) {Q'}^{-1}.
\end{align*}
We show $Q=Q'$ only, for $P=P'$ follows from this by using
\eqref{eq:gamma} and Definition \ref{def:pre-cry} (iii).
{}From \eqref{eq:GiGjGi=GjGiGj} it suffices to check
\begin{align*}
&\vphi_{i}(\boldsymbol{x}) \vphi_{j}(e^c_i(
\boldsymbol{x})) = \vphi_{i}(e^d_j(
\boldsymbol{x})) \vphi_{j}(e^{cd}_i e^d_j(
\boldsymbol{x})), \\
&\vphi_{j}(e^c_i(
\boldsymbol{x})) \vphi_{i}(e^{cd}_j e^c_i(
\boldsymbol{x})) =
\vphi_{j}(\boldsymbol{x}) \vphi_{i}(e^d_j(
\boldsymbol{x})), \\
&\frac{c^{-1}-1}{\vphi_{i}(\boldsymbol{x})} +
\frac{d^{-1}-1}{\vphi_{i}(e^{cd}_j e^c_i(
\boldsymbol{x}))} =
\frac{(cd)^{-1}-1}{\vphi_{i}(e^d_j(
\boldsymbol{x}))}.
\end{align*}
Let us prove the first identity.
Using Lemma \ref{lem:eps-phi-decomp} (b) for the multiple product case we 
have
\begin{align*}
\vphi_{j}(e^{cd}_i e^d_j(
\boldsymbol{x})) &=
\frac{
\vphi_{ij}(e^{cd}_i e^d_j(\boldsymbol{x})) +
\vphi_{ji}(e^{cd}_i e^d_j(\boldsymbol{x}))}
{\vphi_i(e^{cd}_ie^d_j(\boldsymbol{x}))}\\
&=
\frac{
cd \vphi_{ij}(\boldsymbol{x}) +
d \vphi_{ji}(\boldsymbol{x})}
{cd \vphi_{i}( e^d_j(\boldsymbol{x}))} .
\end{align*}
Therefore
\begin{align*}
\vphi_{i}( e^d_j(\boldsymbol{x}))
\vphi_{j}(e^{cd}_i e^d_j(
\boldsymbol{x})) &=
c^{-1} (c \vphi_{ij}(\boldsymbol{x}) + \vphi_{ji}(\boldsymbol{x})) \\
& =
c^{-1} \vphi_{i}(e^{c}_i (\boldsymbol{x}))
\vphi_{j}(e^{c}_i (\boldsymbol{x})) \\
& =
\vphi_{i}(\boldsymbol{x})
\vphi_{j}(e^{c}_i (\boldsymbol{x})).
\end{align*}

The proof of the second identity is similar.
Let us consider the third one.
By similar calculation we have
\[
\vphi_{i}(e^{d}_j (
\boldsymbol{x})) =\frac{
\vphi_{ij}(\boldsymbol{x}) +
d \vphi_{ji}(\boldsymbol{x})}
{d \vphi_{j}(\boldsymbol{x})},\quad
\vphi_{i}(e^{cd}_j e^c_i(
\boldsymbol{x})) =
\frac{
\vphi_{ij}(\boldsymbol{x}) +
d \vphi_{ji}(\boldsymbol{x})}{d}
\times
\frac{
c\vphi_i(\boldsymbol{x})}
{c\vphi_{ij}(\boldsymbol{x}) +
\vphi_{ji}(\boldsymbol{x})}.
\]
Therefore the identity to be proved is written as
\begin{displaymath}
\frac{c^{-1}-1}{\vphi_{i}(\boldsymbol{x})} +
\frac{(d^{-1}-1)d}{
\vphi_{ij}(\boldsymbol{x}) +
d \vphi_{ji}(\boldsymbol{x})} \cdot
\frac{
c \vphi_{ij}(\boldsymbol{x}) +
\vphi_{ji}(\boldsymbol{x})}
{c \vphi_{i}(\boldsymbol{x})} =
\frac{((cd)^{-1}-1)d \vphi_{j}(\boldsymbol{x})}{
\vphi_{ij}(\boldsymbol{x}) +
d \vphi_{ji}(\boldsymbol{x})},
\end{displaymath}
or equivalently
\begin{displaymath}
(c^{-1}-1) (\vphi_{ij}(\boldsymbol{x}) +
d \vphi_{ji}(\boldsymbol{x})) +
(1-d) (\vphi_{ij}(\boldsymbol{x}) +
c^{-1} \vphi_{ji}(\boldsymbol{x})) =
(c^{-1}-d ) \vphi_i(\boldsymbol{x}) \vphi_j(\boldsymbol{x}).
\end{displaymath}
This equality holds since $\vphi_i(\boldsymbol{x}) \vphi_j(\boldsymbol{x}) 
=
\vphi_{ij}(\boldsymbol{x}) + \vphi_{ji}(\boldsymbol{x})$.

The proof of \eqref{eq:ijeqji} is much simpler due to
Lemma \ref{lem:eps-phi-decomp} (a) for the multiple product case.
The proof is completed.
\end{proof}
We present an alternative proof of this theorem
in appendix \ref{sec:appA}.
\section{Tropical $R$}\label{sec:R}

\subsection{Definition and basic properties}
\label{subsec:properties}
Let $M(x,z)=A(x)+zB(x)+z^2C(x)$ be the $2n$ by $2n$ matrix defined in
the previous section.
Consider the relation
\begin{equation}
M(x,z)M(y,z)=M(x',z)M(y',z).
\label{eq:mmeqmmx}
\end{equation}
Given $x, y$, we regard it as a system of simultaneous equations for 
the $4n-2$ unknowns $x', y'$ 
containing a generic parameter $z$.
\begin{theorem}
\label{pr:mainmain}
There is a unique solution of (\ref{eq:mmeqmmx})
under the constraints $\ell(x)=\ell(y') \ne \ell(y)=\ell(x')$.
The map $(x,y) \rightarrow (x',y')$ is a birational transformation.
\end{theorem}

\begin{definition}\label{def:RRRR}
The unique solution $(x',y')$ of 
(\ref{eq:mmeqmmx}) specified in Theorem 
\ref{pr:mainmain} is denoted by $R(x,y)$, and the 
birational transformation $(x,y) \rightarrow R(x,y)$ is called 
the {\em tropical $R$} for the geometric $D^{(1)}_n$-crystal $\B$.
\end{definition}

In Theorem \ref{pr:mainmain} uniqueness is immediate 
{}from Theorem \ref{th:uniqueness}. 
To show the existence, we introduce an explicit 
birational transformation ${\tilde R}$ 
in section \ref{subsec:tilde_R}.
We then prove in section \ref{subsec:tilde_R_equals_R} 
that ${\tilde R}$ actually solves (\ref{eq:mmeqmmx}), hence 
$R = {\tilde R}$.

Here we exhibit basic properties of the
tropical $R$ by assuming its existence.
Let $x=(x_1,\ldots,\xb_1)$ and
$y=(y_1,\ldots,\yb_1)$ be sets of variables for the geometric $D^{(1)}_n$-crystal
$\B$.
On the pair $(x,y)$ we introduce mutually commuting 
involutive automorphisms:
\begin{equation}\label{eq:pairauto}
(x,y)^{\sigma_1} = (x^{\sigma_1},y^{\sigma_1}), \quad
(x,y)^{\sigma_n} = (x^{\sigma_n},y^{\sigma_n}), \quad
(x,y)^{*} = (x^{*},y^{*})
\end{equation}
in terms of $\sigma_1, \sigma_n$ specified in 
(\ref{eq:sigma1o}) and (\ref{eq:sigmano}).
The automorphism $\ast$ is the only one that 
mixes $x$ and $y$ and defined as
\begin{align}
\label{eq:ast}
\ast &: x_i \longleftrightarrow \yb_i, \quad
\xb_i \longleftrightarrow y_i \quad (1 \leq i \leq n-1),
\quad x_n \longleftrightarrow y_n.
\end{align}
\begin{lemma}
\label{lem:JMJequalsMxxx}
Set 
\begin{equation*}
J_* = E_{n,n} + E_{n+1,n+1} + \sum_{i=1}^{n-1} 
(E_{i,2n+1-i}+E_{2n+1-i,i}).
\end{equation*}
Then we have
\begin{align*}
&J_* \,^t\!M(y,z) J_* = M(x^*,z),\;J_* \,^t\!M(x,z) J_* = M(y^*,z),\\
&J_*\, {\mathcal D}_1(y) J_* = {\mathcal D}_2(x^\ast), \; 
J_*\, {\mathcal D}_1(x) J_* = {\mathcal D}_2(y^\ast).
\end{align*}
\end{lemma}
\begin{proof}
By noting that  $x^* = \tau(y^{\sigma_n})$,
$y^* = \tau(x^{\sigma_n})$, and $J_* = J J_n$, 
the first two relations follow from 
(\ref{eq:JMJequalsMno2}) and (\ref{eq:JMJequalsMno3}).
The latter relations are derived by considering 
the 1st column of the former.
\end{proof}

Let us present basic properties of $R$.
\begin{prop}
\label{lem:sigmaonexx}
\begin{align}
\label{eq:sigma_1_commute_with_R}
R((x,y)^{\sigma_1}) &=(R(x,y))^{\sigma_1}, \\
\label{eq:sigma_n_commute_with_R}
R((x,y)^{\sigma_n}) &=(R(x,y))^{\sigma_n}, \\
\label{eq:ast_commute_with_R}
R((x,y)^*) & =(R(x,y))^*.
\end{align}
\end{prop}
\begin{proof}
Write the left (resp.~right) side 
of (\ref{eq:sigma_1_commute_with_R}) as
$((x^{\sigma_1})',(y^{\sigma_1})')$
(resp.~$((x')^{\sigma_1},(y')^{\sigma_1})$).
Then we have
\begin{align*}
M((x^{\sigma_1})',z) M((y^{\sigma_1})',z) &=
M(x^{\sigma_1},z) M(y^{\sigma_1},z) \\
&= J_1(z) M(x,z) M(y,z) J_1(z) \\
&= J_1(z) M(x',z) M(y',z) J_1(z) \\
&= M((x')^{\sigma_1},z) M((y')^{\sigma_1},z).
\end{align*}
Here we have used Lemma \ref{lem:JMJequalsM}.
Due to Theorem \ref{th:uniqueness} we have 
(\ref{eq:sigma_1_commute_with_R}).
(\ref{eq:sigma_n_commute_with_R}) can be shown similarly.
Write the left (resp.~right) side of
(\ref{eq:ast_commute_with_R}) as
$((x^{*})',(y^{*})')$
(resp.~$((x')^{*},(y')^{*})$).
Then we have
\begin{align*}
M((x^{*})',z) M((y^{*})',z) &=
M(x^{*},z) M(y^{*},z) \\
&= J_* \,^t\!M(y,z) \,^t\!M(x,z) J_* \\
& = J_* \,^t\!(M(x,z) M(y,z)) J_* \\
& = J_* \,^t\!(M(x',z) M(y',z)) J_* \\
&= J_* \,^t\!M(y',z) \,^t\!M(x',z) J_* \\
&= M((x')^{*},z) M((y')^{*},z).
\end{align*}
Here we have used Lemma \ref{lem:JMJequalsMxxx}.
Again Theorem \ref{th:uniqueness} implies 
(\ref{eq:ast_commute_with_R}).
\end{proof}
\begin{prop}
\label{th:conservesignatures}
$\veps_i(R(x,y)) = \veps_i(x,y), \quad \vphi_i(R(x,y)) = \vphi_i(x,y)$.
\end{prop}
\begin{proof}
For $1 \leq i \leq n-1$, Lemma \ref{lem:epsandphiforaaa} states that
\begin{align*}
\veps_i(x,y) &= \frac{(M(x,0)M(y,0))_{i+1,i}}{(M(x,0)M(y,0))_{i,i}}, \\
\vphi_i(x,y) &= \frac{(M(x,0)M(y,0))_{i+1,i}}{(M(x,0)M(y,0))_{i+1,i+1}},
\end{align*}
hence the claim follows.
The case $i=0$ is reduced to $i=1$ as
\begin{align*}
\veps_0(R(x,y)) &= \veps_0 ((R(x^{\sigma_1},y^{\sigma_1}))^{\sigma_1}) \\
&= \veps_1(R(x^{\sigma_1},y^{\sigma_1})) \\
&= \veps_1(x^{\sigma_1},y^{\sigma_1}) \\
&= \veps_0(x,y).
\end{align*}
Here we have used Lemma \ref{lem:vepsvphisigma} and
Proposition \ref{lem:sigmaonexx}.
Similarly $i=n$ case follows from  $i=n-1$ case.
\end{proof}
\begin{prop} \label{th:eRequalsRe}
\begin{displaymath}
e^c_i R = R e^c_i.
\end{displaymath}
\end{prop}
\begin{proof}
For $e_i^c(x,y)=(e_i^{c_1}(x),e_i^{c_2}(y))$ we write 
\begin{displaymath}
R(e_i^c(x,y))=((e_i^{c_1}(x))',(e_i^{c_2}(y))'),
\end{displaymath}
and for $R(x,y) = (x',y')$,
\begin{displaymath}
e_i^c(R(x,y))=(e_i^{c_1'}(x'),e_i^{c_2'}(y')),
\end{displaymath}
where
\begin{math}
c_1'=\frac{c\vphi_i(x')+\veps_i(y')}{\vphi_i(x')+\veps_i(y')},\;
c_2'=\frac{\vphi_i(x')+\veps_i(y')}{\vphi_i(x')+c^{-1}\veps_i(y')}.
\end{math}
By using Theorem \ref{th:MMMequalGMMMG} and Proposition 
\ref{th:conservesignatures}
we have
\begin{align*}
&M((e_i^{c_1}(x))',z)M((e_i^{c_2}(y))',z) =
M(e_i^{c_1}(x),z)M(e_i^{c_2}(y),z) \\
&=G_i\left( \frac{c-1}{z^{\delta_{i0}}\veps_i(x,y)} \right)
M(x,z)M(y,z) \;
G_i\left( \frac{c^{-1}-1}{z^{\delta_{i0}}\vphi_i(x,y)} \right)\\
&=G_i\left( \frac{c-1}{z^{\delta_{i0}}\veps_i(x',y')} \right)
M(x',z)M(y',z) \;
G_i\left( \frac{c^{-1}-1}{z^{\delta_{i0}}\vphi_i(x',y')} \right)\\
&=M(e_i^{c_1'}(x'),z)M(e_i^{c_2'}(y'),z).
\end{align*}
Then Theorem \ref{th:uniqueness} tells
$((e_i^{c_1}(x))',(e_i^{c_2}(y))') = (e_i^{c_1'}(x'),e_i^{c_2'}(y'))$.
\end{proof}

\begin{prop} \label{th:RRequalsId}
$R(R(x,y))=(x,y)$.
\end{prop}
\begin{proof}
Writing $R(x,y)=(x',y')$ and $R(x',y')=(x'',y'')$, we get
$M(x,z)M(y,z) = M(x',z)M(y',z)=M(x'',z)M(y'',z)$ by definition.
Then $(x,y)=(x'',y'')$ owing to Theorem \ref{th:uniqueness}.
\end{proof}

Finally we prove the Yang-Baxter equation.
Let $\B\times\B\times\B$ be the 3-fold product of the 
geometric $D^{(1)}_n$-crystal $\B$.
Denote by $R_{12}$ (resp.~$R_{23}$)
the birational map on $\B \times \B \times \B$
that acts on the first (resp. last) two components as $R$ and the other
single component trivially.
\begin{prop}
\label{th:yang-baxter}
We have $R_{12} R_{23} R_{12} = R_{23} R_{12} R_{23}$
as a birational map on $\B \times \B \times \B$.
\end{prop}
\begin{proof}
Let $(x^1,x^2,x^3) \in \B \times \B \times \B$.
Writing $(u^{1},u^{2},u^{3}) = R_{12} R_{23} R_{12} (x^1,x^2,x^3)$ 
and $(v^{1},v^{2},v^{3}) = R_{23} R_{12} R_{23}(x^1,x^2,x^3)$, 
one has
\begin{align*}
M(u^{1},z) M(u^{2},z) M(u^{3},z) &=
M(x^1,z) M(x^2,z) M(x^3,z) \\
&=M(v^{1},z) M(v^{2},z) M(v^{3},z).
\end{align*}
Thus $(u^{1},u^{2},u^{3}) = (v^{1},v^{2},v^{3})$
follows {}from Theorem \ref{th:uniqueness}.
\end{proof}
\subsection{Birational map $\tilde{R}$}
\label{subsec:tilde_R}
\noindent
Let $V_0=(A(x)A(y))_{2n,1}$. Explicitly we have
\begin{align}
\label{eq:16}
V_0 &= \ell (x) \frac{y_1}{\yb_1} + \ell (x) \sum_{m=2}^{n-1}
\left( \prod_{i=1}^{m-1} \frac{y_i}{x_i} \right)
\left( 1+\frac{y_{m}}{\yb_{m}} \right) +
\left( \prod_{i=1}^{n-1} \xb_i y_i \right) x_n y_n \\
\nonumber
&+ \ell (y) \frac{\xb_1}{x_1} + \ell (y) \sum_{m=2}^{n-1}
\left( \prod_{i=1}^{m-1} \frac{\xb_i}{\yb_i} \right)
\left( 1+\frac{\xb_{m}}{x_{m}} \right) +
\left( \prod_{i=1}^{n-1} \xb_i y_i \right).
\end{align}
Define $V_i \, (1 \leq i \leq n-1)$ by
\begin{align}
&V_i = \frac{\yb_i}{\xb_i}V_{i-1} + (\ell (x) - \ell (y))\left( 
1+\frac{\yb_i}{x_i} \right),
\quad (1 \leq i \leq n-2)
\label{eq:vind}\\
&V_{n-1} = \frac{\yb_{n-1}}{\xb_{n-1}}V_{n-2} +
(\ell (x) - \ell (y)) \left( \frac{1}{y_n } +
\frac{\yb_{n-1}}{ x_{n-1} } \right).
\label{eq:unv}
\end{align}
Define $W_i \, (1 \leq i \leq n-1)$ by
\begin{eqnarray}
&& W_i = V_i V_i^* + (\ell (y) - \ell (x))V_i^* + (\ell (x) - \ell 
(y))V_i,
\quad (1 \leq i \leq n-2)
\label{eq:27c}\\
&& W_{n-1} = V_{n-1} V_{n-1}^{*}.
\label{eq:27b}
\end{eqnarray}
Here $V_i^\ast = V_i((x,y)^\ast)$, and similar notations 
will be used for $\sigma_1, \sigma_n$ from now on.

\begin{definition}
\label{th:tropicalR}
A rational transformation specified below is 
denoted by $\tilde{R}(x,y) = (x',y')$:
\begin{eqnarray*}
&& x_1' = y_1 \frac{V_0^{\sigma_1}}{V_1}, \quad
\xb_1' = \yb_1 \frac{V_0}{V_1}, \\
&& x_i' = y_i \frac{V_{i-1}W_i}{V_iW_{i-1}}, \quad
\xb_i' = \yb_i \frac{V_{i-1}}{V_i}, \quad (2 \leq i \leq n-1)\\
&& x_n' = y_n \frac{V_{n-1}}{V_{n-1}^{*}},  \\
&& y_1' = x_1 \frac{V_0}{V_1^*}, \quad
\yb_1' = \xb_1 \frac{V_0^{\sigma_1}}{V_1^*}, \\
&& y_i' = x_i \frac{V_{i-1}^*}{V_i^*}, \quad
\yb_i' = \xb_i \frac{V_{i-1}^* W_i}{V_i^* W_{i-1}}, \quad (2 \leq i \leq 
n-1)\\
&& y_n' = x_n \frac{V_{n-1}^*}{V_{n-1}}.
\end{eqnarray*}
\end{definition}

To provide an explicit formula for $V_i$'s, we introduce 
\begin{align*}
\theta_{i,j}(x,y) &=
\begin{cases}
\ell(x) {\displaystyle
\prod_{k=j+1}^{i} \frac{\yb_k}{\xb_k}} & \mbox{for} \quad
1 \leq j \leq i, \\
\ell(y) {\displaystyle
\prod_{k=i+1}^{j} \frac{\xb_k}{\yb_k}} & \mbox{for} \quad
i+1 \leq j \leq n-2,
\end{cases}
\\
\theta'_{i,j}(x,y) &= \ell(x) \left( \prod_{k=1}^{i} \frac{\yb_k}{\xb_k}
\right) \left(
\prod_{k=1}^{j} \frac{y_k}{x_k} \right) \quad
\mbox{for} \quad j=1,\ldots,n-2,\\
\end{align*}
\begin{align*}
\eta_{i,j}(x,y) &=
\begin{cases}
\ell(x) {\displaystyle 
\left( \prod_{k=j+1}^{i} \frac{\yb_k}{\xb_k} \right) 
\left( \frac{\yb_{j}}{x_{j}} \right)} &\mbox{for} \quad 1 \leq j \leq i,\\
\ell(y) {\displaystyle 
\left( \prod_{k=i+1}^{j} \frac{\xb_k}{\yb_k} \right) 
\left( \frac{\yb_{j}}{x_{j}} \right)} &\mbox{for} \quad i+1 \leq j \leq 
n-1,\\
\ell(y) {\displaystyle \left( \prod_{k=i+1}^{n-1} \frac{\xb_k}{\yb_k} 
\right)
x_n }&\mbox{for} \quad j=n,
\end{cases}
\\
\eta'_{i,j}(x,y) &=
\begin{cases}
\ell(x) {\displaystyle \left( \prod_{k=1}^{i} \frac{\yb_k}{\xb_k}
\right) \left(
\prod_{k=1}^{j} \frac{y_k}{x_k} \right) \left( \frac{x_{j}}{\yb_{j}} \right) 
}
\quad \mbox{for} \quad 1 \leq j \leq n-1,\\
\ell(x) {\displaystyle \left( \frac{\ell(x)}{\ell(y)} \right)
^{\delta_{i,n-1}}
\left( \prod_{k=1}^{i} \frac{\yb_k}{\xb_k}
\right) \left(
\prod_{k=1}^{n-1} \frac{y_k}{x_k} \right) \left( \frac{1}{x_{n}} \right) }
\quad \mbox{for} \quad j=n.
\end{cases}
\end{align*}
\begin{prop}
\label{th:subtractionfree}
For $0 \leq i \leq n-1$ one has
\begin{displaymath}
V_i = \sum_{j=1}^{n-2} \left( \theta_{i,j}(x,y) + \theta'_{i,j}(x,y) 
\right)
+ \sum_{j=1}^{n} \left( \eta_{i,j}(x,y) + \eta'_{i,j}(x,y) \right).
\end{displaymath}
\end{prop}
\begin{proof}
For $i=0$ the above expression agrees with (\ref{eq:16}).
Note that $\eta_{0,n}(x,y)=
\left( \prod_{i=1}^{n-1} \xb_i y_i \right) x_n y_n$ and
$\eta'_{0,n}(x,y)=
\left( \prod_{i=1}^{n-1} \xb_i y_i \right)$.
For $1 \leq i \leq n-2$ the recursion relation 
(\ref{eq:vind}) is valid since the nonzero 
contribution to the difference $V_i - V_{i-1}\yb_i/\xb_i$ only comes 
{}from $\theta_{i,i}(x,y)-\frac{\yb_i}{\xb_i}\theta_{i-1,i}(x,y) =
\ell(x)-\ell(y)$ and
$\eta_{i,i}(x,y)-\frac{\yb_i}{\xb_i}\eta_{i-1,i}(x,y) =
(\ell(x)-\ell(y))\frac{\yb_i}{x_i}$.  
Similarly (\ref{eq:unv}) is checked by 
$\eta_{n-1,n-1}(x,y)-\frac{\yb_{n-1}}{\xb_{n-1}}\eta_{n-2,n-1}(x,y) =
(\ell(x)-\ell(y))\frac{\yb_{n-1}}{x_{n-1}}$ and
$\eta'_{n-1,n}(x,y)-\frac{\yb_{n-1}}{\xb_{n-1}}\eta'_{n-2,n}(x,y) =
(\ell(x)-\ell(y))\frac{1}{y_{n}}$.
\end{proof}
Proposition \ref{th:subtractionfree} with $i=n-1$ reads
\begin{align}\label{eq:Vnminusoneexplicit}
V_{n-1}
&= x_1 y_1 x_n \left( \prod_{m=2}^{n-1} x_m \yb_m \right)
+ \ell (y) \sum_{m=2}^{n-1}
\left( \prod_{i=m}^{n} \frac{x_i}{y_i} \right)
\left( 1+\frac{y_{m}}{\yb_{m}} \right) +
\ell (y) x_n \\
&+ \xb_1 \yb_1 x_n \left( \prod_{m=2}^{n-1} x_m \yb_m \right)
+ \ell (x) \sum_{m=2}^{n-1}
\left( \prod_{i=m}^{n-1} \frac{\yb_i}{\xb_i} \right)
\left( 1+\frac{\xb_{m}}{x_{m}} \right) +
\ell (x) \frac{1}{y_n} ,
\nonumber
\end{align}
{}from which 
$V_{n-1} = V_{n-1}^{\sigma_1} = (V_{n-1}^*)^{\sigma_n}$ is easily seen.
Similarly from (\ref{eq:16}) one finds 
$V_0 = V_0^{\sigma_n} = V^\ast_0$.
Starting from these properties one can use 
(\ref{eq:vind})--(\ref{eq:27b}) to figure out 
the transformation property of 
$V_i, W_i$  under the commuting automorphisms 
$\sigma_1, \sigma_n$ and $\ast$. 
The result is summarized in

\begin{table}[h]
	\caption{Transformation by automorphisms and $\tilde{R}$.}
	\begin{tabular}[h]{|c|c|c|c|c|}
		\hline
		& $V_0$ & $V_i \, (1 \leq i \leq n-2)$ & $V_{n-1}$ & 
		$W_i \, (1 \leq i \leq n-1)$ \\
		\hline
	$\sigma_1$& $V_0^{\sigma_1}$ & $V_i$ & $V_{n-1}$ & $W_i$ \\
	$\sigma_n$& $V_0$ & $V_i$ & $V_{n-1}^*$ & $W_i$ \\
	$\ast$ & $V_0$ & $V_i^*$ & $V_{n-1}^*$ & 
	$W_i$ \\
	$\tilde{R}$ & $V_0$ & $W_i/V_i^*$& $V_{n-1}$& $W_i$ \\
		\hline
	\end{tabular}
	\label{tab:one}
\end{table}

\noindent
The transformation properties under $\tilde{R}$ will be shown in
Lemma \ref{lem:invariances}.

\begin{prop}
\label{lem:rstar}
\begin{align}
\label{eq:commute_with_sigma_one}
\tilde{R}((x,y)^{\sigma_1})&=(\tilde{R}(x,y))^{\sigma_1},\\
\label{eq:commute_with_sigma_n}
\tilde{R}((x,y)^{\sigma_n}) &= (\tilde{R}(x,y))^{\sigma_n},\\
\label{eq:commute_with_ast}
\tilde{R}((x,y)^*)&=(\tilde{R}(x,y))^*.
\end{align}
\end{prop}
\begin{proof}
Apply Table 1 to Definition \ref{th:tropicalR}.
\end{proof}
\begin{lemma}
\label{lem:factorW1}
$W_1 = V_0 V_0^{\sigma_1}$.
\end{lemma}
\begin{proof}
Consider the identities
\begin{math}
V_0 = \frac{\xb_1}{\yb_1} V_1 + 
(\ell (y) - \ell (x))  \frac{\xb_1}{\yb_1} 
\left( 1 +  \frac{\yb_1}{x_1} \right),
\end{math}
\begin{math}
V_0^{\sigma_1}  =\frac{\yb_1}{\xb_1} V_1^* + (\ell (x) - \ell (y)) 
 \frac{\yb_1}{\xb_1}
\left( 1 +  \frac{\xb_1}{y_1} \right),
\end{math}
which are equivalent to (\ref{eq:vind}) and
$* \circ \sigma_1$ of (\ref{eq:vind}) with $i=1$.
Thus we obtain
\begin{math}
V_0 V_0^{\sigma_1} - V_1 V_1^* -
(\ell (y) - \ell (x)) V_1^* - (\ell (x) - 
\ell (y)) V_1
= (\ell (y) - \ell (x)) \frac{\yb_1}{x_1}V_1^*
         + (\ell (x) - \ell (y)) \frac{\xb_1}{y_1}V_1-
         (\ell (x) - \ell (y))^2 \left( 1 +  \frac{\yb_1}{x_1} \right)
         \left( 1 +  \frac{\xb_1}{y_1} \right).
\end{math}
The right hand side vanishes due to (\ref{eq:vind}) and
$*$ of (\ref{eq:vind}) for $i=1$.
\end{proof}
\begin{lemma}
\begin{equation}
W_i = \frac{x_i \yb_i}{\xb_i y_i}
V_{i-1} V_{i-1}^* + (\ell (x) - \ell (y)) \frac{\yb_i}{y_i}
V_{i-1}^* + (\ell (y) - \ell (x)) \frac{x_i}{\xb_i} V_{i-1},
\label{eq:40b}
\end{equation}
for $1 \leq i \leq n-1$.
\end{lemma}
\begin{proof}
For $i=1$ this relation reduces to
\begin{math}
V_0^{\sigma_1} = \frac{x_1 \yb_1}{\xb_1 y_1} +
(\ell (x) - \ell (y)) \frac{\yb_1}{y_1} +
(\ell (y) - \ell (x)) \frac{x_1}{\xb_1},
\end{math}
which follows from (\ref{eq:vind}) and $\sigma_1$ of
(\ref{eq:vind}) with $i=1$.
For $2 \leq i \leq n-2$ the relation is shown by
substituting (\ref{eq:vind}) and $*$ of (\ref{eq:vind}) into (\ref{eq:27c}). 
Next consider $i=n-1$ case.
{}From (\ref{eq:unv}) and $*$ of (\ref{eq:unv}) one has
\begin{math}
\mbox{LHS} - \mbox{RHS}
= (\ell (x) - \ell (y)) \frac{x_{n-1}  }{y_{n-1} y_n} V_{n-2}^* +
  (\ell (y) - \ell (x)) \frac{\yb_{n-1}  }{\xb_{n-1} x_n} V_{n-2}
 - (\ell (x) - \ell (y))^2
  \left( 1 + \frac{x_{n-1}  }{\yb_{n-1} y_n} \right)
  \left( 1 + \frac{\yb_{n-1}  }{x_{n-1} x_n} \right).
\end{math}
The right hand side vanishes because of
(\ref{eq:unv}) and $* \circ \sigma_n$ of (\ref{eq:unv}).
\end{proof}
\begin{lemma}
\label{lem:W_subtraction_free}
\begin{equation}
\left( \frac{1}{x_i} + \frac{1}{\yb_i} \right) W_i =
\frac{1}{y_i} V_i V_{i-1}^* + \frac{1}{\xb_i} V_{i-1} V_i^*
\quad (1 \leq i \leq n-2).
\label{eq:44b}
\end{equation}
\end{lemma}
\begin{proof}
For $i=1$ this relation reduces to
\begin{math}
\left( \frac{1}{x_1} + \frac{1}{\yb_1} \right)V_0^{\sigma_1}=
\frac{V_1^*}{\xb_1} +\frac{V_1}{y_1},
\end{math}
which follows from $\sigma_1$ of (\ref{eq:vind})
and $\sigma_1 \circ *$ of (\ref{eq:vind}) with $i=1$.
To see the other cases,
substitute (\ref{eq:vind}) and $*$ of (\ref{eq:vind}) into the right
hand side of (\ref{eq:44b}). 
The result is the right hand side of 
(\ref{eq:40b}) multiplied by 
$\left( \frac{1}{x_i} + \frac{1}{\yb_i} \right)$.
\end{proof}
\begin{remark}\label{rem:subtractionfree}
We call a rational function {\it subtraction-free},
if its denominator and numerator are polynomials with 
nonnegative coefficients. 
{}From Proposition \ref{th:subtractionfree}, 
Lemma \ref{lem:factorW1}, \ref{lem:W_subtraction_free}
and (\ref{eq:27b}),  all the functions $V_i, W_i$ 
appearing in Definition \ref{th:tropicalR} 
are subtraction-free Laurent polynomials in $x,y$. 
Thus ${\tilde R}$ is also subtraction-free.
This property will be used in section \ref{subsec:piecewise}.
\end{remark}
%
\begin{prop}
\label{th:main}
$\tilde{R}(x,y) =(x',y')$ solves
the equation $A(x)A(y)=A(x')A(y')$.
\end{prop}
\noindent
A proof by a direct but lengthy calculation is available 
in appendix \ref{sec:proofofaaeqaa}.
\begin{lemma}
\label{lem:invariances}
$V_0, V_0^{\sigma_1}, V_{n-1}, V_{n-1}^{*},
W_i \quad (1 \leq i \leq n-1)$ are invariant under $\tilde{R}$.
$\tilde{R}$ acts as $V_i \rightarrow W_i/V_i^*$ and
$V_i^* \rightarrow W_i/V_i \quad (1 \leq i \leq n-1)$.
\end{lemma}
\begin{proof}
Proposition \ref{th:main} tells that 
$V_0 = (A(x)A(y))_{2n,1}$ is invariant.
Then Proposition \ref{lem:rstar} tells $V_0^{\sigma_1}$
and thereby $W_1 = V_0 V_0^{\sigma_1}$ are also invariant.
Consider $V_1$ and $V_1^*$.
By applying $\tilde{R}(x,y)=(x',y')$
to (\ref{eq:vind}) with $i=1$ we have
\begin{math}
\tilde{R}(V_1) =\frac{\yb_1'}{\xb_1'}V_0 + (\ell (x') - \ell (y'))
\left( 1+\frac{\yb_1'}{x_1'} \right)
=\frac{V_0^{\sigma_1}}{V_1^*} \left[
\frac{\xb_1}{\yb_1 } V_1+ (\ell (y) -\ell (x))
\frac{\xb_1}{V_0^{\sigma_1}}
\left( \frac{V_1^*}{\xb_1} +\frac{V_1}{y_1} \right) \right].
\end{math}
By substituting (\ref{eq:44b}) and (\ref{eq:vind})
with $i=1$ into its right hand side
we find $\tilde{R}(V_1)=V_0^{\sigma_1} V_0/V_1^* = W_1/V_1^*$.
Hence by Proposition \ref{lem:rstar} 
we also have $\tilde{R}(V_1^*)= W_1/V_1$.

The claim for $V_i$ and $V_i^*$ with $2 \leq i \leq n-2$
is proved by induction on $i$.
Suppose that the claim is true for $V_{i-1}$.
By applying $\tilde{R}(x,y)=(x',y')$
to (\ref{eq:vind}) we have
\begin{math}
\tilde{R}(V_i) =
\frac{\yb_i'}{\xb_i'}\frac{W_{i-1}}{V_{i-1}^*} +
(\ell (x') - \ell (y'))
\left( 1+\frac{\yb_i'}{x_i'} \right) 
=\frac{W_i}{V_{i-1} V_i^*} \left[
\frac{\xb_i}{\yb_i } V_i+ (\ell (y) - \ell (x))
\frac{\xb_i}{W_i}
\left( \frac{V_{i-1} V_i^*}{\xb_i} 
\right.\right.
\end{math}
\par\noindent
\begin{math}
\left. \left.
+\frac{V_{i-1}^* V_i}{y_i} \right) 
\right] .
\end{math}
By substituting (\ref{eq:44b}) and (\ref{eq:vind})
into the right hand side
we obtain $\tilde{R}(V_i) = W_i/V_i^*$.
Hence by Proposition \ref{lem:rstar} we also
have $\tilde{R}(V_i^*)= W_i/V_i$ for $2 \leq i \leq n-2$.

Let us consider $W_i$ for $2 \leq i \leq n-1$.
By applying $\tilde{R}(x,y)=(x',y')$ to (\ref{eq:40b}) we have
\begin{math}
\tilde{R}(W_i) = \frac{y_i \xb_i}{\yb_i x_i}
\frac{(W_i)^2}{V_{i-1} V_{i-1}^*} + (\ell (y) - \ell (x)) 
\frac{\xb_i}{x_i}
\frac{W_i}{V_{i-1}} + (\ell (x) - \ell (y)) \frac{y_i}{\yb_i}
\frac{W_i}{V_{i-1}^*}.
\end{math}
Dividing both hand sides by $W_i$ and substituting
(\ref{eq:40b}) into the right hand side we find
$\tilde{R}(W_i)/W_i = 1$ for $2 \leq i \leq n-1$.

{}From (\ref{eq:unv}) and its $\sigma_n \circ *$, 
the identity
\begin{math}
V_{n-1} \left( \frac{1}{x_{n-1}x_n} + \frac{1}{\yb_{n-1}} \right)
=\frac{V_{n-2}}{\xb_{n-1}} +\frac{V_{n-2}^*}{y_{n-1}y_n}
\end{math} can be derived.
Applying $\tilde{R}$ to this and using the identity
once again we get $\tilde{R}(V_{n-1})=V_{n-1}$.
The $\tilde{R}$ also leaves
$V_{n-1}^*$ invariant since
$\tilde{R}$ and $*$ are commutative.
\end{proof}
\begin{remark}
{}From the explicit formulas (\ref{eq:16}) and
(\ref{eq:Vnminusoneexplicit}) we find that
$V_0^{\sigma_1} = (A(x)C(y))_{1,2n} + (B(x)B(y))_{1,2n}
+(C(x)A(y))_{1,2n}$,
$V_{n-1} = (A(x)B(y))_{n,n+1}+(B(x)A(y))_{n,n+1}$, and
$V_{n-1}^* = (A(x)B(y))_{n+1,n}+(B(x)A(y))_{n+1,n}$.
Therefore they are invariant under the tropical $R$, and thus 
should be so also under ${\tilde R}$.
\end{remark}
\begin{prop}
\label{th:inversion}
$\tilde{R}(\tilde{R}(x,y)) = (x,y).$
\end{prop}
\begin{proof}
Let us illustrate the $x_2$ case.  The other cases are similar.
Under $\tilde{R}(x,y)=(x',y')$ 
$x_2$ becomes $x_2' = y_2 V_1 W_2/(V_2 W_1)$.
By applying it once again it changes into
\begin{math}
y_2' \frac{R(V_1)R(W_2)}{R(V_2)R(W_1)} =
x_2 \frac{V_1^*}{V_2^*} \frac{(W_1/V_1^*)W_2}{(W_2/V_2^*)W_1} =
x_2.
\end{math}
\end{proof}
Proposition \ref{th:inversion} is the inversion relation of  ${\tilde R}$, 
and will play a role in the proof of Lemma \ref{th:afterward}.

Set
\begin{equation}
\tilde{V}_{j} =
\frac{1}{V_{j}^*}
\left( W_{j} + \frac{\xb_{j+1}}{x_{j+1}}\, W_{j+1} \right), \quad
(1 \leq j \leq n-2).
\label{eq:67apr}
\end{equation}
%
%
\begin{lemma}
\label{lem:26apr}
\begin{eqnarray*}
\tilde{V}_1 &=& V_0 \frac{\yb_1}{\xb_1} \overline{Y}_2
+ (\ell (x) - \ell (y)) \left[
\frac{\yb_2}{y_2} \overline{X}_2 +
\frac{\yb_1}{x_1} \overline{Y}_2 \right],\\
\tilde{V}_j &=& \tilde{V}_{j-1} \frac{\yb_j}{\xb_j}
\frac{\overline{Y}_{j+1}}{\overline{Y}_j} +
(\ell (x) - \ell (y)) \left[
\frac{\yb_{j+1}}{y_{j+1}} \overline{X}_{j+1} +
\frac{\yb_j}{\xb_j}
\frac{\overline{X}_j}{\overline{Y}_j} 
\overline{Y}_{j+1} \right],
\\
&& \qquad \qquad \qquad \qquad (2 \leq j \leq n-2)
\nonumber
\end{eqnarray*}
where $\overline{X}_i = 1+\xb_i/x_i$, 
$\overline{Y}_i = 1+\yb_i/y_i$ for $1 \leq i \leq n-1$.
\end{lemma}
\begin{proof}
{}From (\ref{eq:27c}) and (\ref{eq:40b}) we get
\begin{displaymath}
\tilde{V}_j =
\overline{Y}_{j+1} V_{j} + (\ell (x) - \ell (y))
\left( \frac{\xb_{j+1} \yb_{j+1}}{x_{j+1} y_{j+1}} -1 \right).
\end{displaymath}
This allows us to write $V_j$ in terms of $\tilde{V}_j$.
Substitution of them into (\ref{eq:vind}) leads to the desired 
relations.
\end{proof}

\begin{remark}
\label{rem:tildeV}
The equation 
(\ref{eq:67apr}) gives an alternative definition of $W_j$'s in terms of
$\tilde{V}_j$'s.
Given $W_1$ and $\tilde{V}_{j}$'s one can determine all the other 
$W_j$'s from (\ref{eq:67apr}).
This fact will be used in 
the proof of Lemma \ref{th:yprimesolve}.
\end{remark}
\subsection{Proof of Theorem \ref{pr:mainmain}}
\label{subsec:tilde_R_equals_R}
\noindent
As mentioned after Definition \ref{def:RRRR}, 
Theorem \ref{pr:mainmain} is established 
at the same time with 
the explicit form $R = {\tilde R}$ once it is shown that 
$(x',y') = {\tilde R}(x,y)$ 
solves the defining equation (\ref{eq:mmeqmmx}) 
of the tropical $R$.
In view of $B(x)C(y)+C(x)B(y)=C(x)C(y)=O$, we have the expansion
\begin{align*}
M(x,z)M(y,z) &=A(x)A(y) + z (A(x)B(y)+B(x)A(y)) \\
&+ z^2 (A(x)C(y)+B(x)B(y)+C(x)A(y)).
\end{align*}
Since the equation (\ref{eq:mmeqmmx}) is quadratic with respect to $z$,
it suffices to check it for three distinct values of $z$.
We have already done it at  $z=0$ in Proposition \ref{th:main}.
In what follows we treat the cases 
$z=\ell(x)^{-1}$ and $z=\ell(y)^{-1}$.


Let us put $z= \ell(x)^{-1}=\ell(y')^{-1}$ in (\ref{eq:mmeqmmx}).
 Proposition \ref{prop:meqdpd} leads to
\begin{equation}
P {\mathcal D}_1 (x) M(y, \ell(x)^{-1}) {\mathcal D}_1 (y')^{-1}=
{\mathcal D}_2 (x)^{-1} M(x', \ell(x)^{-1}) {\mathcal D}_2 (y') P.
\label{eq:pdmdeqdmdp}
\end{equation}
Thus the both sides is a scalar multiple of $P$.
Denoting the scalar by $\alpha$ we see that 
(\ref{eq:pdmdeqdmdp}) is equivalent to the simultaneous equations
\begin{eqnarray}
&&P {\mathcal D}_1 (x) M(y, \ell(x)^{-1}) =
\alpha P {\mathcal D}_1 (y'),
\label{eq:mateq1} \\
&&M(x', \ell(x)^{-1}) {\mathcal D}_2 (y') P =
\alpha {\mathcal D}_2 (x) P.
\label{eq:mateq2}
\end{eqnarray}
Equation (\ref{eq:mateq2})
will be considered later (Lemma \ref{th:afterward}).
\begin{lemma}
\label{pr:August4_2}
The following $y_i'$'s and $\yb_i'$'s solve the equation (\ref{eq:mateq1}).
\begin{eqnarray}
&& y_i' = \frac{Q_i}{Q_{i+1}} \quad (1 \leq i \leq n), \quad
\yb_1' = \ell (x) \frac{Q_{2n}}{Q_{2}}, \\
&& \label{eq:QQoverQQ}
\frac{\yb_i'}{y_i'} = - \frac
{\sum_{k=1}^{i} (-1)^k Q_{k} Q_{2n+1-k}}
{\sum_{k=1}^{i-1} (-1)^k Q_{k} Q_{2n+1-k}}, \quad (2 \leq i \leq n-1)
\end{eqnarray}
where $Q_i = (A(x) M(y,\ell (x)^{-1}))_{2n,i}$.
\end{lemma}
\begin{proof}
All elements in the $i$th column of the LHS (resp.~RHS) 
of (\ref{eq:mateq1}) are equal to
$Q_i$ (resp.~$\alpha A(y')_{2n,i}$).
Considering the first column fixes the scalar as
$\alpha = Q_1/A(y')_{2n,1}=Q_1 / \ell (x)$.
Taking the ratios of the $i$th and the $i+1$th column
for $1 \leq i \leq n$ leads to $y_i' = \frac{Q_i}{Q_{i+1}}$.
Similarly the second column and the last column imply
$\yb_1' = \ell (x) \frac{Q_{2n}}{Q_{2}}$, and 
the $(2n-1)$th column does 
$Q_{2n-1} = \alpha \left( 1+ \frac{\yb_2'}{y_2'} \right)\yb_1'
= \left( 1+ \frac{\yb_2'}{y_2'} \right) \frac{Q_1 Q_{2n}}{Q_2}$.
Therefore we have $1+ \frac{\yb_2'}{y_2'}=
\frac{Q_2 Q_{2n-1}}{Q_1 Q_{2n}}$ or equivalently
(\ref{eq:QQoverQQ}) for $i=2$.
Now let us show
\begin{equation}
\label{eq:QQoverQQx}
1+\frac{\yb_i'}{y_i'} = (-1)^{i+1} \frac
{Q_{i} Q_{2n+1-i}}
{\sum_{k=1}^{i-1} (-1)^k Q_{k} Q_{2n+1-k}},
\end{equation}
or equivalently (\ref{eq:QQoverQQ}) for $3 \leq i \leq n-1$.
This claim is checked by induction on $i$.
We have already shown $i=2$ case.
Assume the identity (\ref{eq:QQoverQQ})
with $i$ replaced by $i-1$.
Then we have
\begin{displaymath}
\frac{Q_{2n+1-i}}{Q_{2n+2-i}} =
\frac{\left( 1+\frac{\yb_i'}{y_i'} \right)\frac{\yb_{i-1}'}{y_{i-1}'}}
{\left( 1+\frac{\yb_{i-1}'}{y_{i-1}'} \right)} y_{i-1}'
= \frac{\left( 1+\frac{\yb_i'}{y_i'} \right)
\sum_{k=1}^{i-1} (-1)^k Q_{k} Q_{2n+1-k}}
{(-1)^{i+1} Q_{i-1} Q_{2n+2-i}} \frac{Q_{i-1}}{Q_{i}},
\end{displaymath}
which is equivalent to (\ref{eq:QQoverQQx}).
\end{proof}
\begin{lemma}
\label{th:yprimesolve}
The $y_i'$'s and $\yb_i'$'s specified in Lemma \ref{pr:August4_2}
coincide with those
given by $\tilde{R}(x,y)=(x',y')$
in Definition \ref{th:tropicalR}.
\end{lemma}
\begin{proof}
Set
\begin{eqnarray*}
&& ^{\triangle}\!V_0 = Q_1, \quad
^{\triangle}\!V_i^* = x_1 \ldots x_i Q_{i+1} \quad (1 \leq i \leq n-1),
\\
&&
^{\triangle}\!V_{n-1} = x_1 \ldots x_{n-1} x_n Q_{n+1}, \quad
^{\triangle}\!V_0^{\sigma_1} = \ell (x) \frac{x_1}{\xb_1} Q_{2n}, \\
&&
^{\triangle}\!W_i = \ell (x) \frac{x_1 \ldots x_i}{\xb_1 \ldots \xb_i}
(-1)^i {\sum_{k=1}^{i} (-1)^k Q_{k} Q_{2n+1-k}}
\quad (1 \leq i \leq n-1).
\end{eqnarray*}
We show that these quantities coincide with the
same symbols without $^\triangle$.

First note 
$^{\triangle}\!V_0 =(A(x)A(y))_{2n,1}=V_0$ from (\ref{eq:16}).
It is not difficult to see that $^{\triangle}\!V_0$, 
$^{\triangle}\!V_0^{\sigma_1}$,
$^{\triangle}\!V_i^*$, and $^{\triangle}\!V_{n-1}$
satisfy the same recursion relations
as those without $^\triangle$.
For instance consider
\begin{displaymath}
^{\triangle}\!V_i^* - \frac{x_i}{y_i} \,^{\triangle}\!V_{i-1}^*
= x_1 \ldots x_i \sum_{k=1}^{2n} A(x)_{2n,k}
\left[ M(y,\ell (x)^{-1})_{k,i+1} -
\frac{M(y,\ell (x)^{-1})_{k,i}}{y_i} \right].
\end{displaymath}
By using (\ref{eq:defMno1}) and  $A(x)_{2n,i} = \ell(x)/(x_1\cdots x_{i-1})$, 
this becomes 
\begin{displaymath}
^{\triangle}\!V_i^* = \frac{x_i}{y_i}\,^{\triangle}\!V_{i-1}^* + (\ell (y) - 
\ell (x))\left( 1+\frac{x_i}{\yb_i} \right),
\end{displaymath}
which is $*$ of (\ref{eq:vind}).
The other cases are similar.
Note that $V_0^{\sigma_1}$ is related to $V_1^*$
by $* \circ \sigma_1$ of (\ref{eq:vind}) and $V_{n-1}$ to $V_{n-2}^*$
by $* \circ \sigma_n$ of (\ref{eq:unv}).
By repeated use of the recursion relations,
$^{\triangle}\!V_0^{\sigma_1}$, 
$^{\triangle}\!V_i^*$'s and $^{\triangle}\!V_{n-1}$ are related to
$^{\triangle}\!V_0$ in the same way as those without $^{\triangle}$
are to $V_0$, proving the coincidence.

So far we have considered $^{\triangle}\!V_0$, 
$^{\triangle}\!V_0^{\sigma_1}$,
$^{\triangle}\!V_i^*$, and $^{\triangle}\!V_{n-1}$.
Next we treat $^{\triangle}\!W_i$.
Note $^{\triangle}\!W_1 = \, ^{\triangle}\!V_0
\,^{\triangle}\!V_0^{\sigma_1} = V_0 V_0^{\sigma_1} = W_1$.
To verify $^{\triangle}\!W_i = W_i$ for $i \geq 2$
we define
$^{\triangle}\!\tilde{V}_j=\frac{\ell(x)}{\xb_1 \ldots \xb_j}Q_{2n-j}$
for $1 \leq j \leq n-2$.
Then 
\begin{displaymath}
^{\triangle}\!\tilde{V}_{j} =
\frac{1}{^{\triangle}\!V_{j}^*}
\left( ^{\triangle}\!W_{j} + \frac{\xb_{j+1}}{x_{j+1}}\,
^{\triangle}\!W_{j+1} \right)
\end{displaymath}
holds for $1 \leq j \leq n-2$.
This is (\ref{eq:67apr}).
Since we have already shown $^{\triangle}\!W_1 = W_1$
and $^{\triangle}\!V_j^* = V_j^*$, it suffices to show
$^{\triangle}\!\tilde{V}_j=\tilde{V}_j$ to verify
$^{\triangle}\!W_i = W_i$ for $i \geq 2$.
(See Remark \ref{rem:tildeV}.)
Thus it remains to check that 
$^{\triangle}\!V_0$ and
$^{\triangle}\!\tilde{V}_j$'s satisfy the recursion
relations in Lemma \ref{lem:26apr}.
For instance consider
\begin{displaymath}
^{\triangle}\!\tilde{V}_j -
\,^{\triangle}\!\tilde{V}_{j-1} \frac{\yb_j}{\xb_j}
\frac{\overline{Y}_{i+1}}{\overline{Y}_{i}}
= \frac{\ell(x)}{\xb_1 \ldots \xb_j}
\sum_{k=1}^{2n} A(x)_{2n,k} {\mathcal M}(y,\ell (x)^{-1})_{k,j},
\end{displaymath}
where
\begin{math}
{\mathcal M}(y,z)_{k,j} =
M(y,z)_{k,2n-j} -
\yb_j 
\frac{\overline{Y}_{i+1}}{\overline{Y}_{i}}
M(y,z)_{k,2n+1-j}.
\end{math}
By using (\ref{eq:defMno6}) and 
$A(x)_{2n,2n+1-j} = \xb_1 \ldots \xb_{j-1} 
\left(1+\xb_j/x_j\right)$ we obtain the desired relation.
\end{proof}

Let us proceed to $z= \ell(y)^{-1}$ case.
Setting $z= \ell(y)^{-1}=\ell(x')^{-1}$ in (\ref{eq:mmeqmmx})
leads to
\begin{equation}
{\mathcal D}_2 (x')^{-1} M(x, \ell(y)^{-1}) {\mathcal D}_2 (y) P=
P {\mathcal D}_1 (x') M(y', \ell(y)^{-1}) {\mathcal D}_1 (y)^{-1}.
\label{eq:pdmdeqdmdp2}
\end{equation}
As (\ref{eq:pdmdeqdmdp}), (\ref{eq:pdmdeqdmdp2}) is 
equivalent to the simultaneous equations
\begin{eqnarray}
&&M(x, \ell(y)^{-1}) {\mathcal D}_2 (y) P =
\beta {\mathcal D}_2 (x') P,
\label{eq:mateq3} \\
&&P {\mathcal D}_1 (x') M(y', \ell(y)^{-1}) =
\beta P {\mathcal D}_1 (y)
\label{eq:mateq4}
\end{eqnarray}
for some scalar $\beta$.
Equation (\ref{eq:mateq4})
will be considered later (Lemma \ref{th:afterward}).
\begin{lemma}
\label{th:xprimesolve}
The relation $\tilde{R}(x,y)=(x',y')$ 
solves the equation (\ref{eq:mateq3}) with $\beta = Q^*_1 / \ell (y)$.
\end{lemma}
\begin{proof}
Take the transposition of
(\ref{eq:mateq1}) and
multiply $J_*$'s from the both
sides.  Due to Lemma \ref{lem:JMJequalsMxxx}
the result becomes 
\begin{displaymath}
M(x^*,\ell(x)^{-1}) \mathcal{D}_2 (y^*) P =
\alpha \mathcal{D}_2 ((x')^*) P.
\end{displaymath}
{}From Proposition \ref{lem:rstar} we know $(x')^* = (x^*)'$.
Thus replacing $(x,y)$ with 
$(x,y)^\ast$ in this relation yields (\ref{eq:mateq3}).
\end{proof}
\begin{lemma}
\label{th:afterward}
The relation $\tilde{R}(x,y)=(x',y')$
also solves (\ref{eq:mateq2}) with
$\alpha = V_0/\ell(x)$ and (\ref{eq:mateq4}) with
$\beta = V_0/\ell(y)$.
\end{lemma}
\begin{proof}
Recall $Q_1 =Q^*_1= V_0$ and its invariance under $\tilde{R}$
(Lemma \ref{lem:invariances}).
Due to Proposition \ref{th:inversion} we can 
exchange $(x,y)$ with $(x',y')$ in (\ref{eq:mateq1}) and
(\ref{eq:mateq3}), which become (\ref{eq:mateq4}) and
(\ref{eq:mateq2}), respectively.
Then the claim follows from Lemmas 
\ref{th:yprimesolve} and \ref{th:xprimesolve}.
\end{proof}
\begin{proof}[Proof of Theorem \ref{pr:mainmain}]
Lemmas \ref{pr:August4_2}, \ref{th:yprimesolve}, \ref{th:xprimesolve}, 
and \ref{th:afterward} establish that $(x',y')=\tilde{R}(x,y)$
solves the equation (\ref{eq:mmeqmmx})
at $z = \ell(x)^{-1}$ and $\ell(y)^{-1}$.
With the $z=0$ case verified in Proposition \ref{th:main}, 
we conclude that it solves (\ref{eq:mmeqmmx}) for any $z$ 
as argued in the beginning of the present subsection.
The level constraints $\ell(x)=\ell(y') \ne \ell(y)=\ell(x')$
are obviously satisfied. 
Due to Theorem \ref{th:uniqueness} there is no other
solution.
\end{proof}
\subsection{Properties of $V_i, W_i$ under $e^c_i$}
Here we list the transformation property of 
the functions $V_i, W_i$ under $e^c_i$.
Set
\begin{align*}
&\omega_i = \frac{c \vphi_i(x) + \veps_i(y)}{\vphi_i(x) + \veps_i(y)},
\quad
\psi_i = \frac{c \vphi_i(x) + c \veps_i(y)}{c \vphi_i(x) + \veps_i(y)}
\quad (0 \leq i \leq n),\\
&\Omega_i = 
\begin{cases}
\frac{x_{i+1} + \omega_i \xb_{i+1}}{x_{i+1} + \xb_{i+1}} &
(1 \leq i \leq n-1), \\
\frac{x_2 + \omega_0 \xb_2}{x_2 +\xb_2} &( i=0 ),
\end{cases}
\quad
\Psi_i = 
\begin{cases}
\frac{y_{i+1} + \psi_i \yb_{i+1}}{y_{i+1} + \yb_{i+1}} &
(1 \leq i \leq n-1), \\
\frac{y_2 + \psi_0 \yb_2}{y_2 +\yb_2} &( i=0 ).
\end{cases}
\end{align*}
For a function $F=F((x,y))$
in the variables $(x,y)$
we write 
$R(F) = F(R(x,y))$ and $e^c_i (F) = F(e^c_i(x,y))$.

\begin{prop}\label{pr:V0energylike}
Let $(x',y') = R(x,y)$.
Then we have
\begin{equation*}
V_0(e^c_i(x,y)) = V_0(x,y) 
\left(
\frac{c \vphi_0(x')+ \veps_0(y')}{\vphi_0(x') + \veps_0(y')}
\frac{c \vphi_0(x) + \veps_0(y)}{c \vphi_0(x) + c \veps_0(y)}
\right)^{\delta_{i 0}} .
\end{equation*}
\end{prop}
\begin{proof}
First suppose $1 \le i \le n$.
{}From Theorem \ref{th:MMMequalGMMMG} the relation 
$({\tilde x}, {\tilde y}) = e^c_i(x,y)$  implies 
\begin{math}
M(\tilde{x},z) M(\tilde{y},z) = 
G_i(a)
M(x,z) M(y,z)
G_i(b)
\end{math} for some $a,b$.
But \begin{math}V_0 = (M(x,z) M(y,z))_{2n,1}\end{math} is unchanged 
by the multiplication of the $G_i$'s, proving 
$e^c_iV_0 = V_0$.
To show the $i=0$ case, note 
{}from Lemma \ref{lem:Dynkinrightmost} that 
$e^c_0((x,y)^{\sigma_1}) = (e^c_1(x,y))^{\sigma_1}$.
Thus we have
$e^c_0(V_0^{\sigma_1}) = V_0(e^c_0((x,y)^{\sigma_1})) 
= V_0((e^c_1(x,y))^{\sigma_1}) = (e^c_1V_0)^{\sigma_1}
= V_0^{\sigma_1}$.
Therefore, 
\begin{math}
e^c_0 R \left( \frac{x_1}{\xb_1} \right) =
e^c_0 \left(  \frac{y_1}{\yb_1} \frac{V_0^{\sigma_1}}{V_0} \right) =
\frac{y_1}{\yb_1 \psi_0} \frac{V_0^{\sigma_1}}{e^c_0(V_0)}
\end{math}
and
\begin{math}
R e^c_0 \left( \frac{x_1}{\xb_1} \right) =
R \left( \frac{x_1}{\xb_1 \omega_0} \right) =
\frac{y_1}{\yb_1 R(\omega_0)} \frac{V_0^{\sigma_1}}{V_0}.
\end{math}
Now the assertion follows from 
$e^c_0 R = R e^c_0$ (Proposition \ref{th:eRequalsRe}).
\end{proof}
Proposition \ref{pr:V0energylike} indicates that 
$V_0$ is a tropical analogue of the energy function in 
crystal theory, which will be argued in Remark \ref{rem:energy}.

Similar properties can be derived for the other $V_i, W_i$
with the help of Proposition \ref{th:eRequalsRe} and 
Definition \ref{th:tropicalR}.
Besides {\sc Table} \ref{tab:two}, 
the result is summarized as 
\begin{displaymath}
e^c_i(V_j) = 
\begin{cases}
V_j & (j \ne i), \\
V_i \frac{R(\Omega_i)}{\Psi_i} & ( j = i), 
\end{cases}
\quad
e^c_i(W_j) = 
\begin{cases}
W_j & (j \ne i), \\
W_i \frac{R(\omega_i)}{\psi_i} & ( j = i), 
\end{cases}
\end{displaymath}
for $2 \leq j \leq n-2$ and $0 \leq i \leq n$.
For instance the $j=i$ case of the former (resp.~the latter) relation
can be obtained by applying $e^c_i R = R e^c_i$ on
$\xb_1 \cdots \xb_i$ (resp.~$(\xb_1 \cdots \xb_i)/(x_1 \cdots x_i)$).

\begin{table}[h]
	\caption{Transformation under $e^c_i$.}
	\begin{tabular}[h]{|c|c|c|c|c|c|}
		\hline
		& $V_0$ & $V_0^{\sigma_1}$ & $V_1$ & $V_{n-1}^*$& $V_{n-1}$\\
		\hline
$e^c_0$		& $V_0 \frac{R(\omega_0)}{\psi_0}$ & $V_0^{\sigma_1}$ &
 $V_1\frac{R(\Omega_0)}{\Psi_0}$  & $V_{n-1}^*$ & $V_{n-1}$\\
$e^c_1$		&  $V_0$ &  $V_0^{\sigma_1}\frac{R(\omega_1)}{\psi_1}$& 
 $V_1\frac{R(\Omega_1)}{\Psi_1}$  & $V_{n-1}^*$  & $V_{n-1}$\\
$e^c_i \, (2 \leq i \leq n-2)$		& $V_0$ & $V_0^{\sigma_1}$ & 
$V_1$ & $V_{n-1}^*$& $V_{n-1}$\\
$e^c_{n-1}$		&  $V_0$ &  $V_0^{\sigma_1}$ &  $V_1$ &  
$V_{n-1}^* \frac{R(\omega_{n-1})}{\psi_{n-1}}$& $V_{n-1}$\\
$e^c_n$		& $V_0$ & $V_0^{\sigma_1}$ & $V_1$ & $V_{n-1}^*$& 
$V_{n-1} \frac{R(\omega_n)}{\psi_n}$\\
		\hline
	\end{tabular}
	\label{tab:two}
\end{table}
%
\subsection{Piecewise linear formula for
the combinatorial $R$}
\label{subsec:piecewise}
Suppose a rational function $F=F(x,y;c)$ is subtraction-free 
in the sense of Remark \ref{rem:subtractionfree} 
with respect to the variables $x=(x_1,x_2,\ldots,\xb_1)$,
$y=(y_1,y_2,\cdots,\yb_1)$ and $c$.
By {\it ultradiscretization} of $F$ we mean the expression obtained by 
replacing $+, \times, /$ with $\max, +, -$, respectively in $F$.
(It is called ``tropicalization" in \cite{BK} on the contrary.)
The procedure is well-defined particularly because 
the both sides of $p(q+r) = pq + pr$ have the equal image
$p+\max(q,r) = \max(p+q,p+r)$.
A way to substantiate it is to consider the limit
$\lim_{\epsilon \rightarrow +0} \epsilon
\log F(e^{x/\epsilon},e^{y/\epsilon};e^{c/\epsilon})$
after replacing $x,y,c$ by the real positive arrays 
$e^{x/\epsilon} :=(e^{x_1/\epsilon},e^{x_2/\epsilon},
\ldots,e^{\xb_1/\epsilon})$,
$e^{y/\epsilon} :=(e^{y_1/\epsilon},e^{y_2/\epsilon},
\cdots,e^{\yb_1/\epsilon})$ and $e^{c/\epsilon}$ 
\cite{TTMS}.
Note that ultradiscretized functions are piecewise linear.
In this subsection we use the same symbol to stand for 
the original (tropical) and 
the ultradiscretized (piecewise linear) objects, 
supposing no confusion might arise.
\begin{example}
\label{ex:August6_2}
\begin{displaymath}
F=\frac{c \vphi_0(x)+\veps_0(y)}{\vphi_0(x)+\veps_0(y)} =
\frac{c \xb_1 (1+\xb_2/x_2) + y_1 (1+y_2/\yb_2)}
{\xb_1 (1+\xb_2/x_2) + y_1 (1+y_2/\yb_2)}
\end{displaymath}
is subtraction-free. Its ultradiscretization reads
\[
F = \max \left( c + \xb_1 + (\xb_2 - x_2)_+,
y_1 + (y_2 - \yb_2)_+ \right)
- \max \left(\xb_1 + (\xb_2 - x_2)_+,
y_1 + (y_2 - \yb_2)_+ \right),
\]
where $(x)_+ := \max(x,0)$.
Moreover if $c=1$, 
\[
F = \theta \left( \xb_1 + (\xb_2 - x_2)_+
\geq y_1 + (y_2 - \yb_2)_+ \right).
\]

\end{example}

Let $(x',y') = R(x,y)$ be the ultradiscretization of 
our tropical $R$. Namely, 
\begin{equation}\label{eq:combiR}
\begin{split}
& x_1' = y_1 +V_0^{\sigma_1} -V_1, \quad
\xb_1' = \yb_1 +V_0 -V_1, \\
& x_i' = y_i +V_{i-1}+W_i -V_i -W_{i-1}, \quad
\xb_i' = \yb_i +V_{i-1} -V_i, \quad (2 \leq i \leq n-1)\\
& x_n' = y_n +V_{n-1} -V_{n-1}^{*},  \\
& y_1' = x_1 +V_0 -V_1^*, \quad
\yb_1' = \xb_1 +V_0^{\sigma_1} -V_1^*, \\
& y_i' = x_i +V_{i-1}^* -V_i^*, \quad
\yb_i' = \xb_i +V_{i-1}^* +W_i -V_i^* -W_{i-1}, \quad (2 \leq i \leq 
n-1)\\
& y_n' = x_n +V_{n-1}^* -V_{n-1}.
\end{split}
\end{equation}
Here $V_i, W_i$ are given by 
\begin{align}
\label{eq:piecewise}
&V_i =  \max
\left( \{ \theta_{i,j}(x,y),
\theta'_{i,j}(x,y) | 1 \leq j \leq n-2 \} \right.
\\
\nonumber
& \left. \qquad \qquad \cup
\{ \eta_{i,j}(x,y), \eta'_{i,j}(x,y) | 1 \leq j \leq n \} \right),
\nonumber\\
&W_1 = V_0 + V_0^{\sigma_1},\;\; W_{n-1}=V_{n-1}+V_{n-1}^*, 
\nonumber\\
&W_i=\max \left( V_i+V_{i-1}^*-y_i, V_{i-1}+V_i^*-\xb_i \right)
+ \min (x_i, \yb_i ),\; (2 \leq i \leq n-2),  \nonumber
\end{align}
where 
\begin{align*}
\theta_{i,j}(x,y) &=
\begin{cases}
\ell(x) +{\displaystyle
\sum_{k=j+1}^{i} (\yb_k - \xb_k)} & \mbox{for} \quad
1 \leq j \leq i, \\
\ell(y) +{\displaystyle
\sum_{k=i+1}^{j} (\xb_k - \yb_k)} & \mbox{for} \quad
i+1 \leq j \leq n-2,
\end{cases}
\\
\theta'_{i,j}(x,y) &= \ell(x) + \sum_{k=1}^{i} (\yb_k - \xb_k)+
\sum_{k=1}^{j} (y_k - x_k)  \quad
\mbox{for} \quad j=1,\ldots,n-2,
\end{align*}
\begin{align*}
\eta_{i,j}(x,y) &=
\begin{cases}
\ell(x) +{\displaystyle  \sum_{k=j+1}^{i} (\yb_k -\xb_k)  +
\left( \yb_{j}- x_{j} \right)} &\mbox{for} \quad 1 \leq j \leq i,\\
\ell(y) +{\displaystyle  \sum_{k=i+1}^{j} (\xb_k - \yb_k)  +
\left( \yb_{j}- x_{j} \right)} &\mbox{for} \quad i+1 \leq j \leq n-1,\\
\ell(y) +{\displaystyle  \sum_{k=i+1}^{n-1} (\xb_k - \yb_k) +
x_n }&\mbox{for} \quad j=n,
\end{cases}
\\
\eta'_{i,j}(x,y) &=
\begin{cases}
\ell(x)+ {\displaystyle  \sum_{k=1}^{i} (\yb_k -\xb_k)+
\sum_{k=1}^{j} (y_k -x_k)  + \left( x_{j} -\yb_{j} \right) }
\quad \mbox{for} \quad 1 \leq j \leq n-1,\\
\ell(x)+ {\displaystyle
\delta_{i,n-1} \left( \ell(x) - \ell(y) \right)+
 \sum_{k=1}^{i} (\yb_k -\xb_k)+
\sum_{k=1}^{n-1} (y_k -x_k)  -x_{n} }
\quad \mbox{for} \quad j=n,
\end{cases}
\end{align*}
\begin{equation*}
\ell(x) = \sum_{k=1}^{n} x_k + \sum_{k=1}^{n-1} \xb_k ,\quad
\ell(y) = \sum_{k=1}^{n} y_k + \sum_{k=1}^{n-1} \yb_k .
\end{equation*}
The operations $\sigma_1, \sigma_n$ and $\ast$ in (\ref{eq:combiR}) 
are still defined by (\ref{eq:pairauto}) and (\ref{eq:ast}).

Let us recall the $D^{(1)}_n$-crystal $B_l \; (l \in \Z_{>0})$
{}from \cite{KKM}. It is a finite set 
\[
B_l = \{ (\zeta_1,\ldots,\zeta_n,\overline{\zeta}_n,\ldots,
\overline{\zeta}_1) \in \Z_{\ge 0}^{2n} \mid 
\sum_{i=1}^n(\zeta_i + \overline{\zeta}_i) = l, \; \;
\zeta_n\overline{\zeta}_n = 0 \}
\]
endowed with certain maps ${\tilde e}_i, {\tilde f}_i\, (0 \le i \le n)$
as described in section 5.4 of the mentioned paper.
The set $B_l$ is in one to one correspondence with 
another set
\begin{align*}
B'_l = \{ (x_1,\ldots,x_n,\overline{x}_{n-1},\ldots,
\overline{x}_1)& \in \Z^{2n-1} \mid 
x_i,\xb_i\ge0\text{ for }1\le i\le n-1,\\
&x_n\ge-\min(x_{n-1},\xb_{n-1}),
\sum_{i=1}^{n-1}(x_i + \overline{x}_i) + x_n = l\}
\end{align*}
via $x_i = \zeta_i, \;\;  \overline{x}_i 
= \overline{\zeta}_i \;  (1 \le i \le n-2),\; 
x_{n-1}=\zeta_{n-1}+\overline{\zeta}_{n}, \;
x_{n}=\zeta_{n}-\overline{\zeta}_{n}, \; 
\xb_{n-1}=\overline{\zeta}_{n-1}+\overline{\zeta}_{n},\; 
\zeta_n = \max(0,x_n),\; \overline{\zeta}_n=\max(0,-x_n),\; 
\zeta_{n-1} = x_{n-1} + \min(0,x_n),\; 
\overline{\zeta}_{n-1} = \xb_{n-1} + \min(0,x_n)$.
Note that $x_n$ can be negative. Below we give the explicit crystal structure
of $B'_l$, which is also obtained by ultradiscretizing the geometric crystal
structure of $\B$ in section \ref{subsec:defofB}.
\begin{eqnarray*}
\veps_0(x)&=&x_1+(x_2-\xb_2)_+,\quad
\vphi_0(x)=\xb_1+(\xb_2-x_2)_+,\\
\veps_i(x)&=&\xb_i+(x_{i+1}-\xb_{i+1})_+,\quad
\vphi_i(x)=x_i+(\xb_{i+1}-x_{i+1})_+\quad(i=1,\ldots,n-2),\\
\veps_{n-1}(x)&=&x_n+\xb_{n-1},\quad \vphi_{n-1}(x)=x_{n-1},\\
\veps_n(x)&=&\xb_{n-1},\quad \vphi_n(x)=x_{n-1}+x_n,\\
\et{0}^c(x)&=&(x_1-\xi_2,x_2+\xi_2-c,\ldots,\xb_2+\xi_2,\xb_1-\xi_2+c)  
,\\
\et{i}^c(x)&=&(\ldots,x_i-\xi_{i+1}+c,x_{i+1}+\xi_{i+1}-c,\ldots,
\xb_{i+1}+\xi_{i+1},\xb_i-\xi_{i+1},\ldots)\\
&&\hspace{7.5cm}\quad(i=1,\ldots,n-2),\\
\et{n-1}^c(x)&=&(\ldots,x_{n-1}+c,x_n-c,\ldots),\\
\et{n}^c(x)&=&(\ldots,x_n+c,\xb_{n-1}-c,\ldots),\\
\mbox{where }&&\xi_i=\max(x_i,\xb_i+c)-\max(x_i,\xb_i)\quad(i=1,\ldots,n-1).
\end{eqnarray*}
Here $c$ is an integer. If $c$ is negative, we understand $\et{i}^c$ as 
$\ft{i}^{-c}$. Also, if $\et{i}^c(x)$ does not belong to $B'_l$, one should
assume it to be $0$.

We are interested in the bijection $R:B'_l \otimes B'_k 
\rightarrow B'_k \otimes B'_l$ which commutes with $\et{i}^c$. It is called 
the {\it combinatorial $R$}. See Remark \ref{rem:tensor} for the action of
$\et{i}^c$ on the tensor product. Below we write $(x,y)$ for an element
$x\ot y$ of $B'_l\ot B'_k$. 
We remark that the combinatorial $R: B_l\ot B_k\rightarrow B_k\ot B_l$ has
been given in terms of a certain insertion procedure \cite{HKOT2}.

\begin{theorem}
\label{th:combR}
Restriction of (\ref{eq:combiR}) 
to $x \in B'_{\ell_1}, y \in B'_{\ell_2}$ yields an explicit formula 
$(x,y) \mapsto (x',y')$ 
for the combinatorial 
$R: B'_{\ell_1} \otimes B'_{\ell_2} \rightarrow 
B'_{\ell_2} \otimes B'_{\ell_1}$.
\end{theorem}

\begin{proof}
Set
\[
\tilde{B}'_l=\{(x_1,\ldots,x_n,\overline{x}_{n-1},\ldots,
\overline{x}_1)
\in\Z^{2n-1}\mid
\sum_{i=1}^{n-1}(x_i+\xb_i)+x_n=l\}.
\]
Let $E_i^c$ be the map on $\tilde{B}'_{\ell_1}\ot \tilde{B}'_{\ell_2}$ obtained 
by ultradiscretizing $e_i^c$ on $\B\times\B$, let $R_0$ be the map 
$\tilde{B}'_{\ell_1}\ot\tilde{B}'_{\ell_2}\rightarrow
\tilde{B}'_{\ell_2}\ot\tilde{B}'_{\ell_1}$ given by ultradiscretizing the tropical
$R$, and let $R:B'_{\ell_1}\ot B'_{\ell_2}\rightarrow 
B'_{\ell_2}\ot B'_{\ell_1}$, be the combinatorial $R$. We are to show $R(x,y)=
R_0(x,y)$ if $(x,y)\in B'_{\ell_1}\ot B'_{\ell_2}$. 

First recall that $R$ is uniquely characterized by
\begin{itemize}
\item[(i)] $R((\ell_1,0,\ldots,0),(\ell_2,0,\ldots,0))=
           ((\ell_2,0,\ldots,0),(\ell_1,0,\ldots,0))$,
\item[(ii)] If $(x,y),\et{i}^c(x,y)\in B'_{\ell_1}\ot B'_{\ell_2}$, then
            $R\et{i}^c(x,y)=\et{i}^cR(x,y)$.
\end{itemize}
Let us check (i) for $R_0$. Assume $x=x_0=(\ell_1,0,\ldots,0)$, 
$y=y_0=(\ell_2,0,\ldots,0)$. 
Then we have
$\theta_{i,j} = \ell_1 \, (1 \leq j \leq i), \, = \ell_2 \, (i+1 \leq
j \leq n-2)$;
$\theta'_{i,j}=\ell_2 \, (1 \leq j \leq n-2)$;
$\eta_{i,j} = 0 \, (j=1, i \ne 0), \, = \ell_2 - \ell_1 \,
(j=1, i=0), \, = \ell_1 \, (1 < j \leq i), \, = \ell_2 \,
(i+1 \leq j \leq n)$, and
$\eta'_{i,j} = \ell_1 + \ell_2 \,(j=1), \, = \ell_2 \,
(1 < j \leq n-1 \, \mbox{or} \, j=n, i \ne n-1),
\, =\ell_1 \, (j=n, i=n-1)$.
Therefore we have
$V_0 = V_0^{\sigma_1} = V_i = V_i^* = \ell_1 + \ell_2$ and
$W_i = 2(\ell_1 + \ell_2)$ for $1 \leq i \leq n-1$.
Thus we obtain
$R(x_0,y_0)=((\ell_2 ,0,\ldots,0),(\ell_1 ,0,\ldots,0))$.

For $(x,y)\in B'_{\ell_1}\ot B'_{\ell_2}$ let $m$ be the minimum such that
$\et{i_p}^{c_p}\cdots\et{i_1}^{c_1}(x_0,y_0)\in B'_{\ell_1}\ot B'_{\ell_2}$
($1\le p\le m$), $(x,y)=\et{i_m}^{c_m}\cdots\et{i_1}^{c_1}(x_0,y_0)$
for some integers $c_1,\ldots,c_m$. We prove $R=R_0$ on $B'_{\ell_1}\ot 
B'_{\ell_2}$ by induction on $m$. The $m=0$ case is proven above. Assume 
the $m-1$ case. Note that 
\begin{eqnarray*}
&&E_i^c(x,y)=\et{i}^c(x,y)\text{ if }(x,y),\et{i}^c(x,y)\in 
B'_{\ell_1}\ot B'_{\ell_2},\\
&&E_i^cR_0(x,y)=R_0E_i^c(x,y)\text{ for }(x,y)\in 
\tilde{B}'_{\ell_1}\ot\tilde{B}'_{\ell_2}.
\end{eqnarray*}
The second equality is obtained by ultradiscretizing Proposition 
\ref{th:eRequalsRe}. Thus we get 
\begin{eqnarray*}
R_0(x,y)&=&R_0E_{i_m}^{c_m}\et{i_{m-1}}^{c_{m-1}}\cdots\et{i_1}^{c_1}(x_0,y_0)\\
&=&E_{i_m}^{c_m}R_0\et{i_{m-1}}^{c_{m-1}}\cdots\et{i_1}^{c_1}(x_0,y_0)\\
&=&\et{i_m}^{c_m}R\et{i_{m-1}}^{c_{m-1}}\cdots\et{i_1}^{c_1}(x_0,y_0)\\
&=&R(x,y).
\end{eqnarray*}
The proof is completed.
\end{proof}

\begin{remark}\label{rem:energy}
By the ultradiscretization of 
Proposition \ref{pr:V0energylike}, we have 
\begin{align*}
V_0(\tilde{e}_i(x, y)) &= V_0(x, y) +
\left( \theta (\vphi_0(x') \geq \veps_0(y')) -
\theta(\vphi_0(x) < \veps_0(y)) \right)\delta_{i 0} \\
&= 
\begin{cases}
V_0(x, y) +1&\mbox{ if }i=0,\ \vphi_0(x)\geq\veps_0(y),\ 
\vphi_0(x')\geq\veps_0(y'),\\
V_0(x, y) -1&\mbox{ if }i=0,\ \vphi_0(x)<\veps_0(y),\ 
\vphi_0(x')<\veps_0(y'),\\
V_0(x, y)&\mbox{ otherwise},
\end{cases}
\end{align*}
where $\vphi_0(x) = \xb_1 + (\xb_2 - x_2)_+, \; 
\veps_0(y) = y_1 + (y_2 - \yb_2)_+$.
See Example \ref{ex:August6_2}.
This recursion relation is known to characterize 
an important function on $B'_{\ell(x)} \otimes B'_{\ell(y)}$ called 
the {\em energy function} up to an additive constant.
See [NaY, Definition 3.2].
Thus we conclude that (\ref{eq:16}) is the tropical analogue 
of the energy function. 
In our formula (\ref{eq:piecewise}), it is so normalized that
$\min\{V_0(x,y) \mid (x,y) \in B'_{\ell(x)} \otimes B'_{\ell(y)} \}
= | \ell(x) - \ell(y)|$.
If we put $\ell(x)=\ell(y)=0$, the expression 
(\ref{eq:piecewise}) with $i=0$ reduces to 
the one in \cite{KKM}
associated to the crystal $B_\infty \otimes B_\infty$.
\end{remark}

\appendix
\section{Alternative proof of Theorem \ref{th:ijieqjij}}\label{sec:appA}
\noindent
Let $X=\sum_{1 \leq j \leq i \leq 2n} x_{i,j} E_{i,j}$ be a lower 
triangular
matrix satisfying
\begin{equation}
\label{eq:Liegr}
X S \, ^t \! X S =S \, ^t \! X S X = E,
\end{equation}
where $S$ is the matrix given by (\ref{eq:signmatrix}).
\begin{lemma}
$G_k(u)XG_k(v)^{-1}$ for $1 \leq k \leq n-1$
is lower-triangular if and only if
\begin{align}
\label{eq:uvrel1}
u x_{k+1,k+1} &= v x_{k,k}+u v x_{k+1,k},\\
\label{eq:uvrel2}
u x_{2n+1-k,2n+1-k} &= v x_{2n-k,2n-k}+u v x_{2n+1-k,2n-k}.
\end{align}
\end{lemma}
\begin{proof}
Put $X'=G_k(u)XG_k(v)^{-1}$. The elements of the matrix $X'=(x'_{i,j})$
read
\begin{align*}
x'_{i,j} &= x_{i,j}+u (\delta_{i,k}x_{k+1,j} +
\delta_{i,2n-k}x_{2n+1-k,j})
-v (\delta_{k,j-1}x_{i,k}+\delta_{2n+1-k,j}x_{i,2n-k}) \\
&-uv (
\delta_{i,k}\delta_{j,k+1}x_{k+1,k} +
\delta_{i,2n-k}\delta_{j,k+1}x_{2n+1-k,k} +
\delta_{i,2n-k}\delta_{j,2n+1-k}x_{2n+1-k,2n-k}).
\end{align*}
The only non-zero upper-triangular elements are
\begin{align*}
x'_{k,k+1} &= u x_{k+1,k+1}-v x_{k,k}-u v x_{k+1,k},\\
x'_{2n-k,2n+1-k} &= u x_{2n+1-k,2n+1-k}-v x_{2n-k,2n-k}-u v
x_{2n+1-k,2n-k}.
\end{align*}
Hence the relations (\ref{eq:uvrel1}) and (\ref{eq:uvrel2}) follow.
\end{proof}
Note that the following parameterization solves the
relation (\ref{eq:uvrel1}):
\begin{equation}
\label{eq:uvc}
u=\frac{x_{k,k}}{x_{k+1,k}}(c-1), \quad
v=\frac{x_{k+1,k+1}}{x_{k+1,k}}\frac{(c-1)}{c}.
\end{equation}
It also solves the relation (\ref{eq:uvrel2})
due to (\ref{eq:Liegr}).
In view of this, we define the action $e_i^c \; (1 \leq i \leq n-1)$
on $X$ by
\begin{equation}\label{eq:harry}
e_i^c(X)=G_i\left(\frac{x_{i,i}}{x_{i+1,i}}(c-1)\right)X
G_i\left(\frac{x_{i+1,i+1}}{x_{i+1,i}}\frac{(c-1)}{c}\right)^{-1}.
\end{equation}
Note that we have $e_i^c(X) S \, ^t \! e_i^c(X) S =
S \, ^t \! e_i^c(X) S e_i^c(X) = E$.

\begin{prop}
\label{prop:eee}
\begin{equation*}
e_i^d e_{i+1}^{cd} e_i^c(X)=e_{i+1}^c e_i^{cd} e_{i+1}^d(X)
\quad \hbox{for }\;  1 \leq i \leq n-2.
\end{equation*}
\end{prop}
\begin{proof}
Setting \begin{math}
Y=e_i^c(X), \, Z=e_{i+1}^{cd}(Y), \, Y'=e_{i+1}^d(X), \,
Z'=e_{i}^{cd}(Y'), 
\end{math}
we introduce 
\begin{equation}\label{eq:uv}
\begin{split}
&u_1=\frac{X_{i,i}}{X_{i+1,i}}(c-1), \quad
u_2=\frac{Y_{i+1,i+1}}{Y_{i+2,i+1}}(cd-1), \quad
u_3=\frac{Z_{i,i}}{Z_{i+1,i}}(d-1),\\
&u'_1=\frac{X_{i+1,i+1}}{X_{i+2,i+1}}(d-1), \quad
u'_2=\frac{Y'_{i,i}}{Y'_{i+1,i}}(cd-1), \quad
u'_3=\frac{Z'_{i+1,i+1}}{Z'_{i+2,i+1}}(c-1),\\
&v_1=\frac{X_{i+1,i+1}}{X_{i+1,i}}\frac{(c-1)}{c}, \quad
v_2=\frac{Y_{i+2,i+2}}{Y_{i+2,i+1}}\frac{(cd-1)}{cd}, \quad
v_3=\frac{Z_{i+1,i+1}}{Z_{i+1,i}}\frac{(d-1)}{d},\\
&v'_1=\frac{X_{i+2,i+2}}{X_{i+2,i+1}}\frac{(d-1)}{d}, \quad
v'_2=\frac{Y'_{i+1,i+1}}{Y'_{i+1,i}}\frac{(cd-1)}{cd}, \quad
v'_3=\frac{Z'_{i+2,i+2}}{Z'_{i+2,i+1}}\frac{(c-1)}{c}.
\end{split}
\end{equation}
{}From (\ref{eq:harry}) we know
\begin{align*}
&e_i^d e_{i+1}^{cd} e_i^c(X)=G_i(u_3)G_{i+1}(u_2)G_i(u_1)X
\left(G_i(v_3)G_{i+1}(v_2)G_i(v_1)\right)^{-1},\\
&e_{i+1}^c e_{i}^{cd} e_{i+1}^d(X)=G_{i+1}(u'_3)G_{i}(u'_2)G_{i+1}(u'_1)X
\left(G_{i+1}(v'_3)G_{i}(v'_2)G_{i+1}(v'_1)\right)^{-1}.
\end{align*}
A direct calculation leads to
\begin{equation*}
\begin{array}l
\displaystyle u_1=\frac{x_{i,i}}{x_{i+1,i}}(c-1),\quad
\displaystyle u_2=\frac{x_{i+1,i}x_{i+1,i+1}}{w_1}(cd-1),\quad
\displaystyle u_3=\frac{x_{i,i}}{x_{i+1,i}}\frac{w_1}{w_2}(d-1),\\[4mm]
\displaystyle u'_1=\frac{x_{i+1,i+1}}{x_{i+2,i+1}}(d-1),\quad
\displaystyle u'_2=\frac{x_{i,i}x_{i+2,i+1}}{w_2}(cd-1),\quad
\displaystyle u'_3=\frac{x_{i+1,i+1}}{x_{i+2,i+1}}\frac{w_2}{w_1}(c-1),\\[4mm]
w_1=(1-c) x_{i+2,i}x_{i+1,i+1}+cx_{i+1,i}x_{i+2,i+1}, \\[3mm]
w_2=(d-1)x_{i+2,i}x_{i+1,i+1}+x_{i+1,i}x_{i+2,i+1}.
\end{array}
\end{equation*}
Then the relations
\begin{equation*}
u'_1=\frac{u_2 u_3}{u_1+u_3}, \quad
u'_2=u_1+u_3, \quad
u'_3=\frac{u_1 u_2}{u_1+u_3},
\end{equation*}
can be checked, and from (\ref{eq:GiGjGi=GjGiGj}), 
\begin{equation}\label{eq:GGG=GGG1}
G_i(u_3)G_{i+1}(u_2)G_i(u_1)=G_{i+1}(u'_3)G_{i}(u'_2)G_{i+1}(u'_1).
\end{equation}
Similarly, one can check
\begin{equation}\label{eq:GGG=GGG2}
G_i(v_3)G_{i+1}(v_2)G_i(v_1)=G_{i+1}(v'_3)G_{i}(v'_2)G_{i+1}(v'_1).
\end{equation}
\end{proof}

\begin{lemma}
\label{lem:epsandphiforaaa}
Let $X=\left( X_{i,j} \right)_{1 \leq i,j \leq 2n} = A(x^1)\cdots A(x^L)$.
Then we have
\begin{displaymath}
\frac{X_{i+1,i}}{X_{i,i}} = \veps_i (\boldsymbol{x}), \quad
\frac{X_{i+1,i}}{X_{i+1,i+1}} = \vphi_i (\boldsymbol{x})
\quad \hbox{for } \; 1 \leq i \leq n-1.
\end{displaymath}
\end{lemma}
\begin{proof}
We prove by induction on $L$.
If $L=1$ the claim can be checked directly.
Suppose $L > 1$.
Let $\tilde{X}=\left( \tilde{X}_{i,j}
\right)_{1 \leq i,j \leq 2n} = A(x^1)\cdots A(x^{L-1})$.
Since $\tilde{X}$ and $A(x^L)$ are lower-triangular we have
\begin{align*}
X_{i,i} &= \tilde{X}_{i,i} A(x^L)_{i,i}, \\
X_{i+1,i} &= \tilde{X}_{i+1,i} A(x^L)_{i,i}+
\tilde{X}_{i+1,i+1} A(x^L)_{i+1,i}, \\
X_{i+1,i+1} &= \tilde{X}_{i+1,i+1} A(x^L)_{i+1,i+1}.
\end{align*}
Denote $(x^1,\ldots,x^{L-1})$ by $\tilde{\boldsymbol{x}}$.
Suppose
\begin{displaymath}
\frac{\tilde{X}_{i+1,i}}{\tilde{X}_{i,i}} =
\veps_i (\tilde{\boldsymbol{x}}), \quad
\frac{\tilde{X}_{i+1,i}}{\tilde{X}_{i+1,i+1}} =
\vphi_i (\tilde{\boldsymbol{x}}).
\end{displaymath}
Then we have
\begin{align*}
\frac{X_{i+1,i}}{X_{i,i}} &= \frac{\tilde{X}_{i+1,i}}{\tilde{X}_{i,i}}+
\frac{\tilde{X}_{i+1,i+1}A(x^L)_{i+1,i}}{\tilde{X}_{i,i}A(x^L)_{i,i}}
= \veps_i (\tilde{\boldsymbol{x}})+
\frac{\veps_i (\tilde{\boldsymbol{x}}) \veps_i(x^L)}
{\vphi_i (\tilde{\boldsymbol{x}})}
=\veps_i (\boldsymbol{x}),\\
\frac{X_{i+1,i}}{X_{i+1,i+1}} &=
\frac{A(x^L)_{i+1,i}}{A(x^L)_{i+1,i+1}}+
\frac{\tilde{X}_{i+1,i}A(x^L)_{i,i}}{\tilde{X}_{i+1,i+1}A(x^L)_{i+1,i+1}}
= \vphi_i (x^L)+
\frac{\vphi_i (\tilde{\boldsymbol{x}}) \vphi_i(x^L)}
{\veps_i (x^L)}
=\vphi_i (\boldsymbol{x}).
\end{align*}
\end{proof}

\begin{proof}[Proof of Theorem \ref{th:ijieqjij}]
We concentrate on the nontrivial case (\ref{eq:ijieqjij}) with $i<j$.
Note that $(i,j)=(0,2)$ case reduces to
$(i,j)=(1,2)$ as
\begin{align*}
e^d_0 e^{cd}_2 e^c_0 (\boldsymbol{x}) &=
(e^d_1 e^{cd}_2 e^c_1 (\boldsymbol{x}^{\sigma_1}))^{\sigma_1} \\
&=
(e^c_2 e^{cd}_1 e^d_2 (\boldsymbol{x}^{\sigma_1}))^{\sigma_1} \\
&=e^c_2 e^{cd}_0 e^d_2 (\boldsymbol{x}).
\end{align*}
Similarly  $(i,j)=(n-2,n)$ case reduces to 
$(i,j)=(n-2,n-1)$.
Thus it suffices to prove the relations
$P=P'$ and $Q=Q'$ $(1 \leq i \leq n-2)$ for 
$P,P',Q$ and $Q'$ specified in (\ref{eq:PPQQ1})-(\ref{eq:PPQQ4}).

Let $\boldsymbol{x}=(x^1,\ldots,x^L)$ and $X= A(x^1)\cdots A(x^L)$.
Since $X$ is lower triangular and satisfies (\ref{eq:Liegr}), 
we can apply Lemma \ref{prop:eee} keeping all the notations in its proof.
Then (\ref{eq:GGG=GGG1}) and 
(\ref{eq:GGG=GGG2}) are nothing but $P=P'$ and $Q=Q'$, 
respectively.
Let us show the former for example.
By Lemma \ref{lem:epsandphiforaaa} we have
$u_1 = \frac{c-1}{\veps_{i}(\boldsymbol{x})}$ and 
$u'_1=\frac{d-1}{\veps_{i+1}(\boldsymbol{x})}$.
Next we are to show
$u_2=\frac{cd-1}{\veps_{i+1}(e^c_i(
\boldsymbol{x}))}$ and 
$u'_2=\frac{cd-1}{\veps_{i}(e^d_{i+1}(
\boldsymbol{x}))}$.
We illustrate the way along $u_2$.
Setting  $e^c_i(\boldsymbol{x}) = (y^1,\ldots,y^L)$, we get
\begin{align*}
Y = e^c_i(X) &=
G_i\left( u_1 \right)
X \;
G_i\left( v_1 \right)^{-1} \\
&=
G_i\left( u_1 \right)
A(x^1) \cdots A(x^L) \;
G_i\left( v_1 \right)^{-1} \\
&= A(y^1) \ldots A(y^L),
\end{align*}
where the last equality is due to Theorem \ref{th:MMMequalGMMMG}.
Now Lemma \ref{lem:epsandphiforaaa} tells that 
\begin{math}
\frac{Y_{i+2,i+1}}{Y_{i+1,i+1}} = \veps_{i+1} (e^c_i(\boldsymbol{x})).
\end{math}
The remaining relations are similarly verified.
\end{proof}
\section{Proof of Proposition \ref{th:main}}\label{sec:proofofaaeqaa}
\noindent
Set $I=A(x)A(y)$, which should not be confused with the identity matrix $E$.
\begin{lemma}
\begin{align}
0 &= \sum_{m=i}^n I_{m,i} I_{2n+1-m,i} (-1)^m,
\quad (1 \leq i \leq n)
\label{eq:1c} \\
I_{j,i} &= y_i \left( 1+\frac{x_i}{\yb_i} \right)
\left( 1 + \frac{x_j}{\xb_j} \right) \left( \prod_{m=i+1}^{j-1} x_m 
\right)
+ y_i \left( I_{j,i+1} - \delta_{j,i+1} \right),
\label{eq:reccurent1} \\
& \qquad (1 \leq i \leq j-1 \leq n-2)
\nonumber\\
I_{n,n-1} &= \frac{y_{n-1}}{\yb_{n-1}} x_{n-1} x_n +
y_{n-1} I_{n,n},
\label{eq:4c} \\
I_{n,i} &= y_i \left( 1+\frac{x_i}{\yb_i} \right)
\left( \prod_{m=i+1}^{n} x_m \right)+ y_i I_{n,i+1}, \quad
(1 \leq i \leq n-2)
\label{eq:reccurent2}
\\
I_{j,i} &= y_i \left( 1+\frac{x_i}{\yb_i} \right)
\left( \prod_{m=i+1}^{n} x_m \right)
\left( \prod_{m=2n+1-j}^{n-1} \xb_m \right) + y_i I_{j,i+1},
\label{eq:reccurent3} \\
&\qquad (1 \leq i \leq 2n-j \leq n-2)
\nonumber \\
I_{2n+1-i,i} &= y_i \xb_i 
\left( \frac{1}{x_i} + \frac{1}{\yb_i} \right)
\left( x_i x_n \prod_{m=i+1}^{n-1} x_m \xb_m + \yb_i y_n
\prod_{m=i+1}^{n-1} y_m \yb_m
\right) 
\label{eq:reccurent4} \\
& + y_i \xb_i I_{2n-i,i+1}. \qquad (1 \leq i \leq n-2) \nonumber
\end{align}
\end{lemma}

\begin{proof}
{}From Proposition \ref{pr:MSMS}, one has $I S \, ^t \! I S = E$.
In particular 
$\sum_{m=i}^{2n+1-i} I_{m,i} I_{2n+1-m,i}$ $S_{m,2n+1-m} = 0$,
which is equivalent to (\ref{eq:1c}).

To show  (\ref{eq:reccurent1}), (\ref{eq:reccurent2}) and
(\ref{eq:reccurent3}), we consider
\begin{align*}
I_{j,i} - y_i I_{j,i+1} &= \sum_{k=1}^{2n} A(x)_{j,k}
\left( A(y)_{k,i} - y_i A(y)_{k,i+1} \right) \\
&= A(x)_{j,i} \frac{y_i}{\yb_i} + A(x)_{j,i+1} y_i \\
&= y_i \left( 1+\frac{x_i}{\yb_i} \right) \frac{A(x)_{j,i}}{x_i}
-y_i \delta_{j,i+1},
\end{align*}
where (\ref{eq:defMno1}) is used twice.
Repeated use of (\ref{eq:defMno1}) leads to 
$A(x)_{j,i} = A(x)_{j,1}/(x_{i-1}\cdots x_1)$.
Substituting the explicit form (\ref{eq:firstcolumnofA}) of
$A(x)_{j,1}$ we obtain the desired relations.

Note that $I_{n,n-1}=A(x)_{n,n}A(y)_{n,n-1} + A(x)_{n,n-1}A(y)_{n-1,n-1}
=x_n y_{n-1}y_n + x_{n-1}x_n \frac{y_{n-1}}{\yb_{n-1}}$.
(\ref{eq:4c}) follows from this and $I_{n,n}=x_n y_n$.

Let us consider (\ref{eq:reccurent4}).
{}From (\ref{eq:defMno1}) and (\ref{eq:JMJequalsMno3})
it follows that
\begin{align*}
A(x)_{2n+1-i,k} &= A(x)_{2n-i,k} \xb_i +
\begin{cases}
\xb_i/x_i & \mbox{for} \, k=2n+1-i,\\
\xb_i     & \mbox{for} \, k=2n-i,\\
0         & \mbox{for} \, k \ne 2n+1-i,2n-i,
\end{cases}\\
A(y)_{k,i} &= A(y)_{k,i+1} y_i +
\begin{cases}
y_i/\yb_i & \mbox{for} \, k=i,\\
y_i       & \mbox{for} \, k=i+1,\\
0         & \mbox{for} \, k \ne i,i+1.
\end{cases}
\end{align*}
Substitute these relations into
\begin{math}
I_{2n+1-i,i} = \sum_{k=1}^{2n} A(x)_{2n+1-i,k} A(y)_{k,i}.
\end{math}
The result reads 
\begin{align*}
I_{2n+1-i,i} - y_i \xb_i I_{2n-i,i+1} &=
\xb_i y_i \left( A(x)_{2n-i,i} \frac{1}{\yb_i} + A(x)_{2n-i,i+1} \right) 
\\
&+ \xb_i y_i
\left( A(y)_{2n+1-i,i+1} \frac{1}{x_i} + A(y)_{2n-i,i+1} \right) \\
& = \xb_i y_i \left( \frac{1}{x_i} + \frac{1}{\yb_i} \right)
\left( A(x)_{2n-i,i}  + A(y)_{2n+1-i,i+1}  \right),
\end{align*}
where 
$A(x)_{2n-i,i+1}  = A(x)_{2n-i,i}/x_i$
and $A(y)_{2n-i,i+1} = A(y)_{2n+1-i,i+1}/\yb_i$ are used 
in the last equality.
It remains to show
$A(x)_{2n-i,i} = x_i x_n \prod_{m=i+1}^{n-1} x_m \xb_m$ and
$A(y)_{2n+1-i,i+1} = \yb_i y_n \prod_{m=i+1}^{n-1} y_m \yb_m $.
Repeated use of (\ref{eq:defMno1}) leads to 
$A(x)_{2n-i,i} = A(x)_{2n-i,1}/(x_{i-1}\cdots x_1)$.
The former relation is derived 
by substituting  (\ref{eq:firstcolumnofA}) into $A(x)_{2n-i,1}$ here.
The latter follows from the former by (\ref{eq:JMJequalsMno3}).
\end{proof}
\begin{lemma}
\label{lem:August4_1}
\begin{align}
\label{eq:August4_1_1}
I_{j,i} &= I_{2n+1-i,2n+1-j}^{* \circ \sigma_n}
& (1 \leq i,j \leq n), \\
\label{eq:August4_1_2}
I_{n+1,i} &= I_{n,i}^{\sigma_n} & (1 \leq i \leq n-1).
\end{align}
\end{lemma}
\begin{proof}
{}From the proof of Lemma \ref{lem:JMJequalsMxxx}, we know 
$(x,y) = \tau((x^{\sigma_n},y^{\sigma_n})^\ast)$.
Therefore by using (\ref{eq:JMJequalsMno3}) we have
\begin{math}
A(x)A(y) = 
J\,^t \!A(y^{\ast \circ \sigma_n})
\,^t \!A(x^{\ast \circ \sigma_n}) J =
J\,^t \!(A(x^{\ast \circ \sigma_n} )
A(y^{\ast \circ \sigma_n})) J.
\end{math}
This yields the former relation.
The latter relation follows from
$A(x)A(y) = J_n A(x^{\sigma_n})A(y^{\sigma_n}) J_n$ due to 
(\ref{eq:JMJequalsMno2}).
\end{proof}
The following expressions can be checked directly:
\begin{align}
I_{n+1,n} &= 0,\\
I_{i,i} &= \frac{x_iy_i}{\xb_i\yb_i} \quad ( 1 \leq i \leq n-1),
\quad I_{n,n} = x_n y_n,
\label{eq:diag}
\\
I_{n,n-1} &= x_{n-1}x_n y_{n-1}y_n \left( \frac{1}{x_{n-1} } + 
\frac{1}{\yb_{n-1} y_n } \right),
\label{eq:12a}\\
\label{eq:a11}
I_{n+2,n-1} &= y_{n-1}y_n \xb_{n-1}
\left( \yb_{n-1} + x_{n-1}x_n \right)
\left( \frac{1}{x_{n-1} } + \frac{1}{\yb_{n-1} y_n } \right).
\end{align}
\begin{lemma}
\label{lem:1}
Let $(x',y') = \tilde{R}(x,y)$.
Then we have
\begin{equation*}
1+\frac{x_j'}{\xb_j'}=
\frac{V_{j-1}}{W_{j-1}} \left[
\left( 1+\frac{x_j}{\xb_j} \right) V_{j-1}^* + (\ell (y) - \ell (x))
\left( \frac{x_j y_j}{\xb_j \yb_j} -1 \right) \right],
\end{equation*}
for $2 \leq j \leq n-1$.
\end{lemma}
\begin{proof}
By definition of $\tilde{R}$ we obtain
\begin{equation}
1+\frac{x_j'}{\xb_j'} = 
\frac{1}{W_{j-1}} \left( W_{j-1} + \frac{y_j}{\yb_j} W_j \right).
\end{equation}
In the parentheses of the right side, apply (\ref{eq:27c}) to
(\ref{eq:40b}) to rewrite $W_{j-1}$ and $W_j$.
\end{proof}
\begin{lemma}
\label{lem:2}
Let $(x',y') = \tilde{R}(x,y)$.
Then we have
\begin{equation}
y_i' \left( 1 + \frac{x_i'}{\yb_i'} \right) =
y_i \left( 1 + \frac{x_i}{\yb_i} \right) \frac{W_i}{V_i V_i^*},
\end{equation}
for $1 \leq i \leq n-2$.
\end{lemma}
\begin{proof}
First we suppose $2 \leq i \leq n-2$.
By definition of $\tilde{R}$ we have
\begin{eqnarray}
1 + \frac{x_i'}{\yb_i'} &=&
\frac{y_i}{V_i V_{i-1}^*}
\left( \frac{V_i V_{i-1}^*}{y_i} +
\frac{V_{i-1} V_i^*}{\xb_i} \right) \nonumber\\
&=& \frac{y_i}{V_i V_{i-1}^*}
\left( \frac{1}{x_i} + \frac{1}{\yb_i} \right) W_i,
\end{eqnarray}
where the last step is due to (\ref{eq:44b}).
Multiplying the both sides by 
$y_i' = x_i \frac{V_{i-1}^*}{V_i^*}$ we get the
desired relation.
For  $i=1$, we have
\begin{eqnarray}
1 + \frac{x_1'}{\yb_1'} &=&
\frac{y_1}{V_1}
\left( \frac{V_1}{y_1} +
\frac{V_1^*}{\xb_1} \right) \nonumber\\
&=& \frac{y_1}{V_1}
\left( \frac{1}{x_1} + \frac{1}{\yb_1} \right) V_0^{\sigma_1},
\end{eqnarray}
where the last step is due to 
$\sigma_1$ and $* \circ \sigma_1$ of (\ref{eq:vind}).
Now multiply the both sides by 
$y_1' = x_1 \frac{V_0}{V_1^*}$ and 
simplify the result by Lemma \ref{lem:factorW1}.
\end{proof}
\begin{lemma}
\label{lem:3}
Let $(x',y') = \tilde{R}(x,y)$.
Then we have
\begin{equation}
y_i' = y_i \left[ 1 + \frac{\ell (x) - \ell (y)}{V_i^*} \left( 1 + 
\frac{x_i}{\yb_i} \right)
\right],
\end{equation}
for $1 \leq i \leq n-2$.
\end{lemma}
\begin{proof}
This identity is a consequence of
$y_i' = x_i \frac{V_{i-1}^*}{V_i^*}$ and
$*$ of (\ref{eq:vind}).
\end{proof}
It is straightforward to check 
\begin{lemma}
\label{lem:4}
Let $(x',y') = \tilde{R}(x,y)$.
Then we have
\begin{align*}
\left( \prod_{m=i+1}^{j-1} x'_m \right) &=
\left( \prod_{m=i+1}^{j-1} y_m \right) \frac{V_i W_{j-1}}{V_{j-1}W_i},
\quad (1 \leq i \leq j-1 \leq n-2) \\
%
%
\left( \prod_{m=i+1}^{n} x'_m \right)
\left( \prod_{m=2n+1-j}^{n-1} \xb'_m \right)
&=
\left( \prod_{m=i+1}^{n} y_m \right)
\left( \prod_{m=2n+1-j}^{n-1} \yb_m \right)
\frac{V_i V_{2n-j}}{W_i}.\\
\quad (n+1 \leq j \leq 2n-1)
\end{align*}
\end{lemma}

\begin{lemma}
\label{th:1}
Let $i,j$ be integers such that $1 \leq i \leq j-1 \leq n-2$.
Suppose $I_{j,i+1}$ is invariant under $\tilde{R}$.
Then $I_{j,i}$ is also invariant under $\tilde{R}$.
\end{lemma}
\begin{proof}
Let $(x',y') = \tilde{R}(x,y)$.
By applying Lemmas \ref{lem:1}, \ref{lem:2}, \ref{lem:3} and \ref{lem:4}, 
the relation (\ref{eq:reccurent1}) is written as 
\begin{align}
I_{j,i}(x',y') &= y_i \left( 1+\frac{x_i}{\yb_i} \right)
\frac{1}{V_i^*}
\left[
\left( 1+\frac{x_j}{\xb_j} \right) V_{j-1}^* + (\ell (y) - \ell (x))
\left( I_{j,j} -1 \right) \right]
\left( \prod_{m=i+1}^{j-1} y_m \right) \nonumber\\
& +
y_i \left[ 1 + \frac{\ell (x) - \ell (y)}{V_i^*} \left( 1 + 
\frac{x_i}{\yb_i} \right)
\right]
\left( I_{j,i+1} - \delta_{j,i+1} \right).
\end{align}
Thus it reduces to the following Lemma.
\end{proof}
\begin{lemma}
\begin{align}
& \quad
\left[
\left( 1+\frac{x_j}{\xb_j} \right) V_{j-1}^* + (\ell (y) - \ell (x))
\left( I_{j,j} -1 \right) \right]
\left( \prod_{m=i+1}^{j-1} y_m \right)\nonumber\\
& = \left( \prod_{m=i+1}^{j-1} x_m \right)
\left( 1+\frac{x_j}{\xb_j} \right) V_i^* +
(\ell (y) - \ell (x)) \left( I_{j,i+1} - \delta_{j,i+1} \right),
\label{eq:56b}
\end{align}
for $1 \leq i \leq j-1 \leq n-2$.
\end{lemma}
\begin{proof}
By applying
(\ref{eq:reccurent1}) and $*$ of (\ref{eq:vind})
this is transformed into the same relation with $i$ replaced by $i-1$.
Thus it suffices to check (\ref{eq:56b}) for 
$i=j-1$, which is straightforward.
\end{proof}
\begin{lemma}
\label{th:2}
Let $i$ be an integer such that $1 \leq i \leq n-2$.
Suppose $I_{n,i+1}$ is invariant under $\tilde{R}$.
Then $I_{n,i}$ is also invariant under $\tilde{R}$.
\end{lemma}
\begin{proof}
Let $(x',y') = \tilde{R}(x,y)$.
By applying Lemmas \ref{lem:2}, \ref{lem:3} and \ref{lem:4}, 
the relation (\ref{eq:reccurent2}) is written as
\begin{equation}
I_{n,i}(x',y') = y_i \left( 1+\frac{x_i}{\yb_i} \right)
\frac{V_{n-1}}{V_i^*}
\left( \prod_{m=i+1}^{n} y_m \right)
+
y_i \left[ 1 + \frac{\ell (x) - \ell (y)}{V_i^*} \left( 1 + 
\frac{x_i}{\yb_i} \right)
\right] I_{n,i+1}.
\end{equation}
Thus it reduces to the following Lemma.
\end{proof}
\begin{lemma}
\label{lem:01}
For $1 \leq i \leq n-2$
\begin{equation}
\left( \prod_{m=i+1}^n y_m \right) V_{n-1} =
\left( \prod_{m=i+1}^n x_m \right) V_i^* + 
(\ell (y) - \ell (x)) I_{n,i+1}. 
\label{eq:21a}
\end{equation}
\end{lemma}
\begin{proof}
By using
(\ref{eq:reccurent2}) and $*$ of (\ref{eq:vind})
this is transformed into the same relation with $i$ replaced by $i-1$.
Thus it suffices to check (\ref{eq:21a}) for $i=n-2$, 
which reads
\begin{displaymath}
y_{n-1}y_n V_{n-1} = x_{n-1}x_n V_{n-2}^* + (\ell (y) - \ell (x)) 
I_{n,n-1}.
\end{displaymath}
This follows from (\ref{eq:12a}) and $* \circ \sigma_n$ of (\ref{eq:unv}).
\end{proof}

\begin{lemma}
\label{th:3}
Let $i,j$ be integers such that $1 \leq i \leq 2n-j \leq n-2$.
Suppose $I_{j,i+1}$ is invariant under $\tilde{R}$.
Then $I_{j,i}$ is also invariant under $\tilde{R}$.
\end{lemma}
\begin{proof}
Let $(x',y') = \tilde{R}(x,y)$.
By applying Lemmas \ref{lem:2}, \ref{lem:3} and  \ref{lem:4}, 
the relation (\ref{eq:reccurent3}) is written as
\begin{align}
I_{j,i}(x',y') &= y_i \left( 1+\frac{x_i}{\yb_i} \right)
\frac{V_{2n-j}}{V_i^*}
\left( \prod_{m=i+1}^{n} y_m \right)
\left( \prod_{m=2n+1-j}^{n-1} \yb_m \right)
\nonumber \\
& +
y_i \left[ 1 + \frac{\ell (x) - \ell (y)}{V_i^*} \left( 1 + 
\frac{x_i}{\yb_i} \right)
\right] I_{j,i+1}.
\end{align}
Thus it reduces to the following Lemma.
\end{proof}
\begin{lemma}
\begin{align}
V_{2n-j} \left( \prod_{m=i+1}^{n} y_m \right)
\left( \prod_{m=2n+1-j}^{n-1} \yb_m \right) &=
V_{i}^* \left( \prod_{m=i+1}^{n} x_m \right)
\left( \prod_{m=2n+1-j}^{n-1} \xb_m \right)
\nonumber \\
& + (\ell (y) - \ell (x)) I_{j,i+1},
\label{eq:59}
\end{align}
for $1 \leq i \leq 2n-j \leq n-2$.
\end{lemma}
\begin{proof}
By applying
$*$ of (\ref{eq:vind}) and (\ref{eq:reccurent3})
this is transformed into the same relation with $i$ replaced by $i-1$.
Thus it suffices to check (\ref{eq:59}) for $i=2n-j$, which 
is done in the next Lemma.
\end{proof}
\begin{lemma}
\begin{equation}
V_{i} \left( y_n \prod_{m=i+1}^{n-1} y_m \yb_m \right) =
V_{i}^* \left( x_n \prod_{m=i+1}^{n-1} x_m \xb_m \right)
+ (\ell (y) - \ell (x)) I_{2n-i,i+1},
\label{eq:60}
\end{equation}
for $1 \leq i \leq n-2$.
\end{lemma}
\begin{proof}
By applying (\ref{eq:vind}), $*$ of (\ref{eq:vind}) and 
(\ref{eq:reccurent4})
this is transformed into the same relation with $i$ replaced by $i-1$.
Thus it suffices to check (\ref{eq:60}) for $i=1$.
Consider the trivial identity:
\begin{equation*}
V_0 \left( y_n \prod_{m=1}^{n-1} y_m \yb_m \right) =
V_0 \left( x_n \prod_{m=1}^{n-1} x_m \xb_m \right)
+ (\ell (y) - \ell (x)) I_{2n,1}.
\end{equation*}
The $i=1$ case of (\ref{eq:60}) is shown by 
applying (\ref{eq:vind}), $*$ of (\ref{eq:vind}) and 
(\ref{eq:reccurent4}) to this identity.
\end{proof}
\begin{lemma}
\label{th:4}
Let $i$ be an integer such that $1 \leq i \leq n$.
Then $I_{i,i}$ is invariant under $\tilde{R}$.
\end{lemma}
\begin{proof}
This is clear from (\ref{eq:diag}).
\end{proof}
\begin{lemma}
\label{th:5}
$I_{n,n-1}$ is invariant under $\tilde{R}$.
\end{lemma}
\begin{proof}
Let $(x',y') = \tilde{R}(x,y)$.
{}From (\ref{eq:4c}) and Theorem \ref{th:4} we have
\begin{equation}
I_{n,n-1}(x',y') = \frac{x_{n-1}}{\xb_{n-1}} y_{n-1}y_n
\frac{V_{n-2}}{V_{n-1}^{*}} +
x_{n-1}  \frac{V_{n-2}^*}{V_{n-1}^{*}} I_{n,n}.
\end{equation}
By using $\sigma_n$ and $*$ of (\ref{eq:unv})
we recover the right hand side of (\ref{eq:4c}).
\end{proof}
\begin{proof}[Proof of Proposition \ref{th:main}]
It should be shown that all $I_{j,i}$'s are invariant under $\tilde{R}$.
In view of $I_{n+1,n}=0$, $I_{j,i}=0\, (j < i)$, 
(\ref{eq:August4_1_1}) and Proposition \ref{lem:rstar}, 
it suffices to check the case $i+j \leq 2n+1$.
Suppose $1 \leq j \leq n-1$.
Then the claim follows from Lemmas \ref{th:1} and \ref{th:4}.
Suppose $j=n$.
Then the claim follows from Lemmas \ref{th:2}, \ref{th:4} and 
\ref{th:5}.
Suppose $j=n+1$ and $1 \leq i \leq n-1$.
Then the claim follows from $j=n$ 
due to (\ref{eq:August4_1_2}) 
and Proposition \ref{lem:rstar}.
We have already proved the invariance of
$I_{n-1,n-1}, I_{n,n-1}$ and $I_{n+1,n-1}$.
Because of the identity (\ref{eq:1c}),
$I_{n+2,n-1}$ in (\ref{eq:a11}) is also invariant.
Then the claim for $j=n+2$ follows from Lemma \ref{th:3}.
Similarly the invariance of
$I_{n-a,n-a}, I_{n-a+1,n-a},\cdots,I_{n+a,n-a}
\quad (1 \leq a \leq n-1)$
implies the invariance of $I_{n+a+1,n-a}$ by (\ref{eq:1c}), and
the claim for $j=n+a+1, 1 \leq i \leq n-a$ follows from 
Lemma \ref{th:3}. By induction on $a$ the proof is completed.
\end{proof}


\end{document}